\documentclass[a4paper,reqno]{amsart}
\usepackage{mathrsfs}
\usepackage{comment}
\usepackage{newcent} 
\usepackage[T1]{fontenc}
\usepackage[centertags]{amsmath}
\usepackage{amssymb,amsfonts,amsthm}
\usepackage{newlfont}
\usepackage{fancyhdr,fancyvrb}
\usepackage{graphicx}
\usepackage{pinlabel}
\usepackage{nextpage}
\usepackage[dvipsnames,usenames]{color}
\usepackage{esint}
\usepackage{multirow}
\usepackage{mathtools}
\usepackage{Bon-Dav-Mor-14-private}
\usepackage[colorlinks=true, linkcolor=blue, citecolor=ForestGreen]{hyperref}

\numberwithin{equation}{section}

\theoremstyle{definition}
\begingroup
\newtheorem{defin}{Definition}[section]
\newtheorem{remark}[defin]{Remark}
\endgroup
\theoremstyle{plain}
\begingroup
\newtheorem{theorem}[defin]{Theorem}
\newtheorem{lemma}[defin]{Lemma}
\newtheorem{proposition}[defin]{Proposition}

\endgroup

\title[A perturbed Cahn-Hilliard for LB films]{Analysis of a perturbed Cahn-Hilliard model for Langmuir-Blodgett films}
\author{Marco Bonacini}
\address{Institute for Applied Mathematics, University of Bonn, Endenicher Allee 60, 53115 Bonn, Germany}
\email[M.~Bonacini]{bonacini@iam.uni-bonn.de}
\author{Elisa Davoli}
\address{Department of Mathematical Sciences, University of Vienna, Oskar-Morgenstern-Platz 1, 1090 Vienna, Austria}
\email[E.~Davoli]{elisa.davoli@univie.ac.at}
\author{Marco Morandotti}
\address{Fakult\"at f\"ur Mathematik, Technische Universtit\"at M\"unchen, Boltzmannstrasse, 3, 85748 Garching, Germany}
\email[M.~Morandotti~\myenv]{marco.morandotti@ma.tum.de}

\date{\today}



\makeindex

\begin{document}

\begin{abstract}
An advective Cahn-Hilliard model motivated by thin film formation is studied in this paper.
The one-dimensional evolution equation under consideration includes a transport term, whose presence prevents from identifying a gradient flow structure.
Existence and uniqueness of solutions, together with continuous dependence on the initial data and an energy equality are proved by combining a minimizing movement scheme with a fixed point argument. Finally, it is shown that, when the contribution of the transport term is small, the equation possesses a global attractor and converges, as the transport term tends to zero, to a purely diffusive Cahn-Hilliard equation.
\end{abstract}

\maketitle
\tableofcontents

\section{Introduction} \label{sect:intro}
In this paper we study a perturbation of a Cahn-Hilliard-type equation with an advective term. 
Given $L,T>0$, we consider a one-dimensional evolution equation of the form
\begin{equation}\label{intro0}\tag{EQ$_\beta$}
u_t = \Bigl(-u_{xx} + \frac{\partial W}{\partial s}(x,u) \Bigr)_{xx} - \beta u_x,
\end{equation}
where $u\colon[0,L]\times[0,T]\to\R{}$ is a function which, in the applications, represents the concentration of a phase of fluid or can be a measure of an order parameter, $W\colon[0,L]\times\R{}\to\R{}$ is a (possibly asymmetric) space-dependent double-well potential, and $\beta\geq0$ is the advection velocity. 

In the case $\beta=0$, for the choice of the potential $W$ given by
\begin{equation*}
W(x,s)=W(s) \coloneqq \frac{s^4}4-\frac{s^2}2,
\end{equation*}
\eqref{intro0} reduces to the celebrated Cahn-Hilliard equation
\begin{equation}\label{CHE}
u_t = (-u_{xx} + u^3-u)_{xx}.
\end{equation}
Proposed to describe the evolution of the concentrations of two immiscible fluids, equation \eqref{CHE} has been extensively studied, since the seminal paper \cite{CH58}, from the physical, chemical, and mathematical points of view.
Variations of this model have been proposed to describe different phenomena and modes of phase separation (see, e.g., \cite{NC08} and the references therein).
In particular, introducing a skewness in the wells of $W$ has the effect of selecting a preferred state of the system that minimizes the potential energy. This, in turn, can be further tailored to model a variability of the optimal configuration in space by considering a potential $W$ that depends on the spatial variable $x$.
When the fluids under consideration are subject to stirring, advection becomes significant and it can be incorporated in the equation by introducing the transport term $-\beta u_x$ in \eqref{intro0}, see, e.g., \cite{LON2008}. 
Despite its relevance for applications, the results on equation \eqref{intro0} are mainly numerical (see, e.g., \cite{bertozzi, Frank, Liu}), whereas, to the best of our knowledge, a full analytical treatment is still missing. 

We propose in this paper a first step toward a comprehensive analysis of equation \eqref{intro0} by setting up a functional framework in which we prove existence and uniqueness of solutions, continuous dependence on the initial data, and an energy equality, for every $T>0$ and for every $\beta\geq0$; additionally, for $\beta\leq\beta_0(L)$, where $\beta_0$ is a threshold depending only on the size of the system $L$, we prove the existence of an attractor for \eqref{intro0}, the convergence to the solution of (\blu{EQ$_0$}) for $\beta\to0$, and we quantify the distance between such two solutions, in an appropriate metric, in terms of $\beta$.

\smallskip

Our main motivation for investigating equation \eqref{intro0} is in connection with the mechanism of the Langmuir-Blodgett (LB) transfer. Based on the pioneering results of I.~Langmuir \cite{L17} and K.~B.~Blodgett \cite{B35}, the LB transfer is a way to depose a thin film made of amphiphilic molecules on a solid substrate. The breadth of applications that involve LB films is incredibly wide \cite{O92,roberts}. Many aspects of the physical process of the film formation are interesting for scientists, in particular those related to pattern formation \cite{bottom-up,confinement,science}.

LB films are obtained by pulling a solid substrate out of the so-called Langmuir-Blodgett trough, which contains a monolayer of amphiphilic molecules floating free at the surface of a liquid subphase. Depending on its area density, the monolayer can exist in two different phases, the so-called liquid-expanded (LE) and liquid-condensed (LC) phases. We are interested in the situation in which the density is kept constant during the process and the monolayer is in the LE phase. During the transfer of the amphiphilic molecules onto the substrate, the system can undergo a phase transition into the LC phase, which is energetically favoured due to the interaction with the substrate. However, different scenarios are observed experimentally depending on the transfer velocity, which plays the role of the main control parameter: the transfer of a homogeneous LC phase takes place if the velocity is sufficiently small, whereas, if the substrate is pulled out sufficiently fast, alternating patterns of LE and LC phases emerge \cite{bottom-up}. The use of a prestructured substrate can also give more control over the pattern formation.

In \cite{KGFT12,WG} an evolution equation in the form of a generalized Cahn-Hilliard equation has been proposed to govern the deposition process, starting form a more general model which takes into account also the hydrodynamics of the thin liquid layer \cite{KGFC10}. Denoting by $c=c(z,t)$, $z=(x,y)\in[0,L]\times[0,L]$, $t\in[0,T]$, the concentration of the amphiphiles of a LB film during the transfer, the equation reads
\begin{subequations}\label{intro1a}
\begin{equation} \label{intro1}
c_t = \nabla\cdot\bigl( \nabla(-\Delta c - c + c^3 + \zeta(x)) - vc \bigr),
\end{equation}
where $v=(\beta,0)$ is the transfer velocity, and the profile
\begin{equation*}
\zeta(x)\coloneqq -\frac{\zeta_0}{2}\Bigl(1+\tanh\Bigl(\frac{x-x_s}{l_s}\Bigr)\Bigr)\,,
\end{equation*}
with $\zeta_0>0$ fixed, describes the presence of a meniscus of width $l_s$ located at a point $x_s$. The space-dependent term $\zeta$ introduces a skewness in the two wells of the potential: before the meniscus ($x<x_s$) it has essentially no influence ($\zeta\sim0$) and the two wells $c=-1$ (pure LE phase) and $c=1$ (pure LC phase) are at the same level; after the meniscus, the LC phase is preferred.
In \cite{KGFT12,WG} equation \eqref{intro1} is complemented by the boundary conditions 
\begin{equation}\label{intro2}
c|_{x=0}=c_0, \quad
\frac{\partial c}{\partial x}\bigg|_{x=L}=0, \quad
\frac{\partial^2 c}{\partial x^2}\bigg|_{x=0} = 0, \quad
\frac{\partial^2 c}{\partial x^2}\bigg|_{x=L}=0, \quad
c|_{y=0}=c|_{y=L},
\end{equation}
\end{subequations}
and by an initial condition $c(z,0)=\hat{c}_0(z)$. The parameter $c_0\in\R{}$ represents the (constant) concentration of amphiphiles at the surface of the liquid trough, which is modeled as the interface $\{x=0\}$; by consistency we assume $\hat{c}_0(0)=c_0$. In view of the periodicity conditions in the variable $y$, we will consider in the following a reduced problem in dimension $1$. It is convenient to shift the map $c$ in order to have zero Dirichlet boundary condition at $x=0$: we consider the translation $u(x,t) \coloneqq c(x,t)-c_0$, as well as the translation of the initial datum $u_0(x)\coloneqq \hat{c}_0(x)-c_0$, for $x\in[0,L]$, $t\in [0,T]$. With these positions the evolution equation \eqref{intro1a} 
takes the form
\begin{equation*}
u_t = \bigl(-u_{xx} + (u+c_0)^3 - (u+c_0) + \zeta \bigr)_{xx} - \beta u_x,
\end{equation*}
with
\begin{equation*}
\begin{cases}
u(0,t)=0 & t\in (0,T],\\
u_x(L,t)=0 & t\in (0,T],\\
u_{xx}(0,t)=u_{xx}(L,t)=0 & t\in (0,T],\\
u(x,0)=u_0(x) & x\in[0,L].
\end{cases}
\end{equation*}
This reduces the equation to the form \eqref{intro0}, with the particular choice of the potential
\begin{equation}\label{intro6}
W(x,s) \coloneqq \frac14(s+c_0)^4-\frac12(s+c_0)^2+\zeta(x)s.
\end{equation}
The numerical analysis in \cite{KGFT12} (see also \cite{KT14} for the analysis of the bifurcation structure, and \cite{WG,zhu} for the case of a prestructured substrate) highlights the presence of three different regimes, according to the magnitude of the extraction speed $\beta$: for $\beta$ smaller than a first threshold, the concentration of the amphiphiles moves from the first to the second potential well, corresponding to a transition from the LE phase to the LC phase; for $\beta$ bigger than a second threshold, the concentration remains basically unvaried, so that the film is deposed in a homogeneous LE phase; for intermediate values of $\beta$, the concentration varies forming regular patterns, alternating LE and LC phases.

\smallskip

In this work we begin laying the theoretical foundations for the analysis of \eqref{intro0}. The result of this paper is twofold. 

First, in Theorem~\ref{thm:existence} we prove existence and uniqueness of a weak solution (see Definition~\ref{def:weaksol}) to \eqref{intro0}, satisfying an energy equality. Existence of solutions to advective-diffusion equations with the structure in \eqref{intro0} has been already established by means of Galerkin methods (see, e.g., \cite{bertozzi, LON2008}).
The novelty of our approach is in the fact that it provides a new variational discretization and approximation of solutions to \eqref{intro0}. 
Indeed, the proof strategy for Theorem~\ref{thm:existence} relies on the combination of a minimizing movement scheme with a fixed point argument. To be precise, by decoupling the advective and the diffusive terms in \eqref{intro0} we first show in Section~\ref{sect:MM} existence of solutions for the weak formulation of a Cahn-Hilliard equation with a forcing term. The existence of weak solutions for the advective Cahn-Hilliard equation in \eqref{intro0} is obtained in Section~\ref{sect:fixedpoint} by a fixed point argument. 

Our second contribution is the analysis of the regime ``$\beta$ small''. For $\beta$ smaller than a constant $\beta_0$ dependent only on  $L$, we prove in Theorem~\ref{thm:gb-gen-semiflow} that the equation has a global attractor, in the framework of the classical theory of semigroups (see e.g. \cite{Hale,Lad}).
We point out that related existence results of global attractors have been obtained in \cite{GalGrasselli} for a non-local Cahn-Hilliard equation without advective term, in \cite{FrigeriGrasselli2012, FrigeriGrasselliRocca} for a Cahn-Hilliard equation coupled with a Navier-Stokes system, and in \cite{IuorioMelchionna} for a nonlocal Cahn-Hilliard equation with a reaction term. Attractors for the viscous and standard Cahn-Hilliard equation with non-constant mobility have been analyzed in \cite{GMRS,Schimperna}. The case of a doubly nonlinear Cahn-Hilliard-Gurtin system has been studied in \cite{MiranvilleSchimperna} (see also \cite{akagi} for the case of a doubly nonlinear equation with non-monotone perturbations), whereas that of a logarithmic potential is the subject of \cite{CMZ}. A thorough analysis of nonlocal convective Cahn-Hilliard equations with reaction terms has been undertaken in \cite{DellaPortaGrasselli}.

Finally, in Theorem~\ref{prop:beta} we provide a quantitative estimate of the distance of solutions to \eqref{intro0} from those of (\blu{EQ$_0$}) in terms of the advection velocity $\beta$.
In accordance with the numerical results in \cite{WG}, we conjecture the existence of three analogous regimes for $\eqref{intro0}$. The case $\beta>\beta_0$ is expected to be related to pattern formation, and, together with the higher dimensional setting, will be the subject of a forthcoming paper.

\subsection*{Structure of the paper}
The plan of the paper is the following.
In Section~\ref{sect:setting}, we present the functional setting in which we study the problem and, after introducing a suitable notion of weak solution, we state the main existence result. We also explain the strategy of the proof.
In Section~\ref{sect:MM}, we describe the minimizing movement scheme, exploiting the gradient flow structure that is highlighted in Section~\ref{sect:gradientflow}; here we adapt the general results of \cite{AGS} to the case of time-dependent energies, and, to this end, \cite{A1995} will be a fundamental reference.
In Section~\ref{sect:fixedpoint}, we will apply a fixed point argument to construct the sought weak solution to \eqref{intro1} and complete the proof of Theorem~\ref{thm:existence}. Section~\ref{sect:smallbeta} is devoted to the analysis of the regime where the parameter $\beta$ is small. Finally, we collect in Appendix~\ref{sect:appendix} the statements of some technical lemmas which we use throughout the paper.

\subsection*{Notation}
For a function $u\colon[0,L]\times[0,T]\to\R{}$, $u=u(x,t)$, we will usually denote its partial derivatives by $u'(x,t) \coloneqq \frac{\partial u}{\partial x}(x,t)$ and $\dot{u}(x,t) \coloneqq \frac{\partial u}{\partial t}(x,t)$ . We will also denote by $u(t)$ the function $u(\cdot,t)$ for $t$ fixed.

\section{Functional setting, weak formulation, and existence result} \label{sect:setting}

In this section, we discuss the functional framework to solve equation \eqref{intro0} and we state our existence result. After presenting the precise assumptions on the potential $W$, we introduce the functions spaces in which we set the evolution problem and formulate, in Definition~\ref{def:weaksol}, the notion of weak solution. Then we state Theorem~\ref{thm:existence} about the existence, uniqueness, and continuous dependence from initial data of weak solutions, and the energy equality. An equivalent definition of weak solutions is presented in Proposition~\ref{prop:weaksol}: this will be more convenient for the proof of our main theorem. Finally, in Subsection \ref{subs:strategy} we outline the strategy for the proof of Theorem~\ref{thm:existence}.

We assume that the potential
\begin{subequations}\label{sub1W}
\begin{equation} \label{nd121}
W:[0,L]\times\R{}\to\R{} \quad\text{is of class $C^3$}
\end{equation}
and that it satisfies, for every $x\in [0,L]$ and $s\in \R{}$,
\begin{equation} \label{nd122}
W(x,s) \geq - K_0,
\qquad
\partial_s W(x,s)s\geq W(x,s)-K_1,
\end{equation}
\end{subequations}
for some fixed positive constants $K_0, K_1$. Here and in the following, we denote the partial derivative of $W$ with respect to the variable $s$ by $\partial_s W (x,s) \coloneqq \frac{\partial W}{\partial s}(x,s)$. Notice that the potential \eqref{intro6} satisfies conditions \eqref{sub1W}.

We consider the initial/boundary value problem
\begin{subequations} \label{sub2}
\begin{equation}\label{nd102}
\dot{u} = \bigl(-u'' + \partial_s W(x,u) \bigr)'' - \beta u', \qquad (x,t)\in(0,L)\times(0,T),
\end{equation}
\begin{equation}\label{nd103}
\begin{cases}
u(0,t)=0 & t\in (0,T],\\
u'(L,t)=0 & t\in (0,T],\\
\bigl( -u'' + \partial_sW(x,u) \bigr)|_{x=0} =0 & t\in (0,T],\\
\bigl( -u'' + \partial_sW(x,u) \bigr)'|_{x=L} =0 & t\in (0,T],\\
u(x,0)=u_0(x) & x\in[0,L],
\end{cases}
\end{equation}
\end{subequations}
where $u_0\colon [0,L]\to\R{}$ is a given initial condition with $u_0(0)=0$.

\subsection{Functional setup}
We now introduce the function spaces in which we rigorously set the weak formulation of problem \eqref{sub2}. We set
\begin{equation*}
V\coloneqq \{u\in H^1(0,L): u(0)=0\},
\end{equation*}
endowed with the scalar product and norm 
\begin{equation*}
(u,v)_V \coloneqq \int_0^L u'(x)v'(x)\,\de x,\qquad \norma{u}_V^2=\int_0^L (u'(x))^2\de x,
\end{equation*}
obtained as the restriction to $V$ of the scalar product on $H^1(0,L)$
\begin{equation}\label{nd106}
(u,v)_{H^1(0,L)} \coloneqq u(0)v(0)+\int_0^L u'(x)v'(x)\,\de x
\end{equation}
(which is, in turn, equivalent to the standard scalar product of $H^1(0,L)$). We denote by $V'$ the dual space of $V$ and by $\scp\cdot\cdot_{V',V}$ the duality pairing between $V'$ and $V$, for which we will provide a convenient expression in \eqref{nd108}.
With this notation, we are in a position to give the definition of weak solution to \eqref{sub2}.

\begin{defin}[Weak solution of problem \eqref{sub2}] \label{def:weaksol}
Let $u_0\in V$ be a given initial datum. A function $u$ is a weak solution to \eqref{sub2} in $[0,T]$ corresponding to $u_0$ if $u(x,0)=u_0(x)$ for almost every $x\in (0,L)$, and the following conditions hold:
\begin{itemize}
\item[(ws1)] $u \in H^1(0,T;V')\cap L^2(0,T;H^3(0,L))$;
\item[(ws2)] the function
\begin{equation} \label{nd356}
\mu \coloneqq -u'' + \partial_sW(x,u) \in L^2(0,T;V);
\end{equation}
\item[(ws3)] for every $\psi\in V$ and almost every $t\in(0,T)$ we have
\begin{equation}\label{nd358}
\scp{\dot u(t)}{\psi}_{V',V}= - (\mu(t),\psi)_V -\beta\int_0^L u'(x,t)\psi(x)\,\de x \,;
\end{equation}
\item[(ws4)] for almost every $t\in(0,T)$, $u$ satisfies the boundary conditions
\begin{equation} \label{nd357}
 u(0,t)=0,\quad u'(L,t)=0.
\end{equation}
\end{itemize}
\end{defin}

\begin{remark}\label{rem:mu}
Notice that the higher order boundary conditions in \eqref{nd103} are implicitly contained in the weak formulation: indeed, we are imposing also that the function $\mu(t)$ defined in \eqref{nd356} satisfies the boundary conditions $\mu(0,t)=0$ (which follows from the requirement $\mu(t)\in V$) and $\mu'(L,t)=0$ in a weak sense (as the test functions in \eqref{nd358} do not necessarily vanish at $x=L$).
In other words, for almost every $t\in(0,T)$ the function $\mu(t)$  is a weak solution to
\begin{equation*}
\begin{cases}
\mu''(t) = \dot{u}(t)+\beta u'(t) &\text{in the sense of distributions in $(0,L)$,}\\
\mu(0,t)= \mu'(L,t)=0.
\end{cases}
\end{equation*}
\end{remark}

\begin{remark} \label{rmk:interpolation}
The condition (ws1) in Definition~\ref{def:weaksol} automatically implies that the map $u$ is continuous from $[0,T]$ with values in $V$ (after possibly being redefined on a set of measure zero). This follows by applying \cite[Theorem~8.60]{Leo}.
\end{remark}

The first main result of the paper is the following.

\begin{theorem} \label{thm:existence}
Let $W\colon [0,L]\times\R{}\to\R{}$ satisfy assumptions \eqref{sub1W}, and let $u_0\in V$. Given any $T>0$, there exists a unique weak solution $u$ to the problem \eqref{sub2} in $[0,T]$ corresponding to $u_0$, according to Definition~\ref{def:weaksol}. 
Furthermore, setting
\begin{equation} \label{nd359}
\mathcal{E}(v) \coloneqq \frac12\int_0^L |v'(x)|^2\,\de x + \int_0^L W(x,v(x))\,\de x \qquad\text{for }v\in V,
\end{equation}
the following energy equality holds true for almost every $t\in(0,T)$:
\begin{equation} \label{nd360}
\frac{\de}{\de t}\mathcal{E}(u(t)) + \|\mu(t)\|^2_V = -\beta\int_0^L u'(x,t)\mu(x,t)\,\de x,
\end{equation}
where $\mu$ is defined in \eqref{nd356}.
Finally, the solution $u$ depends continuously on the initial data, in the following sense: given any $M>0$, there exists a constant $C_M$ (depending on $M$ and on $L$, $\beta$, $T$, $W$) such that for every initial data $u_0,\bar u_0\in V$ with $\norma{u_0}_V,\norma{\bar{u}_0}_V\leq M$ the corresponding weak solutions $u,\bar u$ satisfy, for all $t\in[0,T]$,
\begin{equation}\label{S1}
\norma{u(t)-\bar{u}(t)}_V \leq C_M \norma{u_0-\bar{u}_0}_V.
\end{equation}
\end{theorem}
The proof of the theorem is given in Section~\ref{sect:fixedpoint}.

\medskip
As $V$ and $V'$ are Hilbert spaces, it is possible to introduce a  representation for the elements of $V'$ which will be useful later on to rephrase the Definition~\ref{def:weaksol} of weak solution in a way that is more convenient for the proof of Theorem~\ref{thm:existence}. According to the Riesz Representation Theorem, for every $\varphi\in V'$ there exists a unique element $z_\varphi\in V$ such that 
\begin{equation}\label{nd108}
\scp{\varphi}{v}_{V',V}=\int_0^L z_\varphi'(x)v'(x)\,\de x \quad\text{for every }v\in V,
\qquad \norma\varphi_{V'}=\norma{z_\varphi}_V.
\end{equation}
The scalar product in $V'$ can be expressed through the one in $V$ as
\begin{equation}\label{nd110}
(\varphi,\psi)_{V'}=(z_\varphi,z_\psi)_V=\int_0^L z_\varphi'(x)z_\psi'(x)\,\de x \qquad\text{for every }\varphi,\psi\in V'.
\end{equation}
Notice that \eqref{nd108} allows us to characterize $z_\varphi$ as the weak solution to
\begin{equation}\label{nd111}
\begin{cases}
-z_\varphi''=\varphi & \text{in the sense of distributions in $(0,L)$,}\\
z_\varphi(0)= z_\varphi'(L)=0.
\end{cases}
\end{equation}

The standard identification of functions in $L^2(0,L)$ with elements of its dual space $(L^2(0,L))'$, given by the Riesz Representation Theorem, can be formally expressed by the map $Q\colon L^2(0,L) \to (L^2(0,L))'$ defined as
\begin{equation}\label{nd112}
\scp{Qu}{v}_{(L^2)',L^2} \coloneqq (u,v)_{L^2}=\int_0^L u(x)v(x)\,\de x.
\end{equation}
In particular, each function of $L^2(0,L)$ can be seen as a linear map acting on $V$, that is, as an element in $V'$. 
In the following, we will always implicitly identify $u\in L^2(0,L)$ with the corresponding element $Qu\in V'$, without indicating the map $Q$.
Hence, from \eqref{nd112}, for every $u\in L^2(0,L)$ we have
\begin{equation} \label{nd112ter}
\scp{u}{v}_{V',V}=\int_0^L u(x)v(x)\,\de x
\qquad\text{for all }v\in V.
\end{equation}
Notice also that it immediately follows that
\begin{equation} \label{nd112bis}
\norma{u}_{V'} \leq L\norma{u}_{L^2(0,L)}
\qquad\text{for all }u\in L^2(0,L),
\end{equation}
since $\norma{v}_{L^2(0,L)}\leq L\norma{v}_V$ for all $v\in V$. This identification of $L^2(0,L)$ with its dual corresponds to the inclusion diagram
\begin{equation*}
V \subseteq L^2(0,L) \equiv (L^2(0,L))' \subseteq V' \,,
\end{equation*}
where the second inclusion is defined by \eqref{nd112ter}.

The following proposition provides an equivalent formulation of the weak equation \eqref{nd358} in terms of the functions $z$ introduced in \eqref{nd108}.

\begin{proposition}\label{prop:weaksol}
Let $u_0\in V$, and let $u$ satisfy conditions \textup{(ws1)}, \textup{(ws2)} and \textup{(ws4)} in Definition~\ref{def:weaksol}. Then condition \textup{(ws3)} holds if and only if for almost every $(x,t)\in(0,L)\times(0,T)$
\begin{equation} \label{nd300}
-u''(x,t) + \partial_sW(x,u(x,t)) = - z_{\dot{u}(t)}(x) - \beta z_{u'(t)}(x) \,.
\end{equation}
\end{proposition}

\begin{proof}
We have to show the equivalence of \eqref{nd358} and \eqref{nd300}, under the assumptions \eqref{nd356} and \eqref{nd357}.
By taking the scalar product in $V$ of \eqref{nd300} with any test function $\psi\in V$, we have for almost every $t\in [0,T]$
\begin{align*}
(\mu(t),\psi)_V
& = - \bigl( z_{\dot{u}(t)+\beta u'(t)} , \psi \bigr)_V 
\stackrel{\eqref{nd108}}{=}  -\scp{\dot u(t)+\beta u'(t)}{\psi}_{V',V} \\
& \stackrel{\eqref{nd112ter}}{=} -\scp{\dot{u}(t)}{\psi}_{V',V} - \beta\int_0^L u'(x,t)\psi(x)\,\de x\,,
\end{align*}
which is \eqref{nd358}. Conversely, assuming that \eqref{nd358} holds and reversing the previous chain of equalities we find that for almost every $t$
\begin{equation*}
(\mu(t),\psi)_V =  - \bigl( z_{\dot{u}(t)+\beta u'(t)} , \psi \bigr)_V \qquad\text{for all }\psi\in V,
\end{equation*}
which implies \eqref{nd300}.
\end{proof}

\begin{remark} \label{rk:same-space}
The dual space $V'$ can be identified with the space
\begin{equation*}
W \coloneqq \{\varphi\in(H^1(0,L))': \scp\varphi1_{(H^1)',H^1}=0\}
\end{equation*}
of functionals that vanish on constant functions of $H^1(0,L)$. Indeed, the restriction map
\begin{equation*}
i\colon W\to V',\qquad\qquad \varphi\mapsto i(\varphi) \coloneqq \varphi|_{V}
\end{equation*}
is an isometry between the two spaces (considering on $H^1(0,L)$ the norm \eqref{nd106} and on $W$ the corresponding dual norm).
\end{remark}

\subsection{Strategy for the proof of Theorem \ref{thm:existence}}\label{subs:strategy}
Our approach to prove the existence of a weak solution is structured in two steps and consists in the combination of a minimizing movement scheme and of a fixed point argument. To briefly explain the strategy, we consider the function space $Y\coloneqq L^2(0,T;V)$ and we carry out the following program. 
Let $u_0\in V$ be a given initial datum.
\begin{itemize}
\item[\textit{Step 1.}] For every fixed $q\in Y$, we construct in Proposition~\ref{prop:compact} a function $u_q\in Y$ solving
\begin{equation}\label{nd116}
\dot u_q=(-u_q''+ \partial_sW(x,u_q))'' - \beta q'
\end{equation}
in the weak sense, namely, recalling Proposition~\ref{prop:weaksol},
\begin{equation*}
-u_q''+\partial_sW(x,u_q) = -z_{\dot u_q} -\beta z_{q'}\,,
\end{equation*}
and attaining the initial datum $u_q(x,0)=u_0(x)$.
Notice that by construction $u_q(0,t)=0$ for all $t$, since $u_q\in Y$.
The solution $u_{q}$ will be sought for by showing that equation \eqref{nd116}, for fixed $q$, possesses in fact a gradient flow structure (see Section~\ref{sect:gradientflow}), which can be approximated by a minimizing movements scheme. 
In Proposition~\ref{prop:uni} we show that the function $u_q$ obtained in this way is unique.
\smallskip
\item[\textit{Step 2.}] We show in Section~\ref{sect:fixedpoint} that the map $S\colon Y\to Y$, $q\mapsto S(q)\coloneqq u_q$ (where $u_q$ is the function constructed in the previous step), has a fixed point, by applying Shaefer's Fixed Point Theorem.
The fixed point solves
\begin{equation*}
-u_q''+ \partial_sW(x,u_q) = -z_{\dot u_q} -\beta z_{u_q'}\,,
\end{equation*}
which is the weak formulation of \eqref{nd116} for $q=u_q$, by Proposition~\ref{prop:weaksol}.
\end{itemize}

\section{Decoupling of diffusive and advective effects} \label{sect:gradientflow}
In this section we show that for every fixed $q\in L^2(0,T;V)$, equation \eqref{nd116} has a gradient flow structure. For $q\in L^2(0,T;V)$ we define $\Psi_q\colon V'\times[0,T]\to\R{}\cup\{\infty\}$ by
\begin{equation}\label{nd119}
\Psi_q(u,t)\coloneqq
\begin{cases}
\displaystyle \frac12\int_0^L \norm{u'(x)}^2\de x+\int_0^L W(x,u(x))\,\de x+\beta(q'(t),u)_{V'} & \text{if $u\in V$,} \\
\infty & \text{elsewhere in $V'$,}
\end{cases}
\end{equation}
where $W$ is the potential satisfying assumptions \eqref{sub1W}. 
Notice that in \eqref{nd119} the functions $u,q'(t)\in L^2(0,L)$ are implicitly identified with elements in $V'$, according to \eqref{nd112ter}.

In the following proposition we highlight the relation between the evolution equation \eqref{nd116} and the gradient flow of the functional $\Psi_q$ with respect to the metric structure of $V'$, namely
\begin{equation*}
\dot{u}_q\in-\partial_{V'} \Psi_q(u_q,t),
\end{equation*}
where $\partial_{V'}$ denotes the subdifferential of $\Psi_q$ in the space $V'$ with respect to its first variable; that is, $\partial_{V'}\Psi_q(u,t)$ is the set of all elements $g\in V'$ such that
\begin{equation*}
\Psi_q(v,t) \geq \Psi_q(u,t) + ( g,v-u )_{V'} \qquad\text{for all }v\in V'
\end{equation*}
(see, e.g., \cite[Definition~4.48]{FonLeo}).

\begin{proposition} \label{prop:subdiff}
Let $q\in L^2(0,T;V)$ and let $\Psi_q$ be the functional defined in \eqref{nd119}. Assume that $u\in H^3(0,L)$ is such that
\begin{equation*}
u(0)=0, \qquad u'(L)=0,
\end{equation*}
and that the function
\begin{equation*}
\mu \coloneqq -u'' + \partial_s W(\cdot,u(\cdot))
\end{equation*}
belongs to the space $V$. Then for a.e. $t\in (0,T)$ there holds
\begin{equation*}
g\in\partial_{V'}\Psi_q(u,t) \quad\Longrightarrow\quad  \scp{g}{\psi}_{V',V}= (\mu,\psi)_V + \beta\int_0^L q'(x,t)\psi(x)\,\de x \quad\text{for all }\psi\in V\,.
\end{equation*}
\end{proposition}

\begin{proof}
Fix $t\in (0,T)$ such that $q(t)\in V$, and let $u$ be as in the statement of the theorem.
Given $g\in\partial_{V'} \Psi_q(u,t)$, we have that for every $w\in V$ and for every $\eps>0$
\begin{equation*}
\Psi_q(u\pm\eps w,t)-\Psi_q(u,t)\geq\pm\eps(g,w)_{V'}.
\end{equation*}
Dividing by $\eps$ and letting $\eps\to0^+$, we obtain
\begin{equation*}
\int_0^L u'(x)w'(x)\,\de x+\int_0^L \partial_s W(x,u(x))w(x)\,\de x+\beta(q'(t),w)_{V'}=(g,w)_{V'}.
\end{equation*}
In view of Lemma~\ref{lem:mu''} below, we find that if $g\in\partial_{V'}\Psi_q(u,t)$ then
\begin{equation*}
\int_0^L z_{w}'(x)\mu'(x) \,\de x + \beta(q'(t),w)_{V'} = (g,w)_{V'}\,.
\end{equation*}
By \eqref{nd108}--\eqref{nd110} this identity can be rewritten as
\begin{equation*}
(\mu,z_w)_V + \beta\scp{q'(t)}{z_w}_{V',V} = \scp{g}{z_w}_{V',V} \qquad\text{for all }w\in V,
\end{equation*}
and in turn by \eqref{nd112ter} (notice that $q'(t)\in L^2(0,L)$)
\begin{equation} \label{nd131bis}
(\mu,z_w)_V + \beta\int_0^L q'(x,t)z_w(x)\,\de x = \scp{g}{z_w}_{V',V} \qquad\text{for all }w\in V.
\end{equation}

To conclude the proof it is sufficient to show that \eqref{nd131bis} holds if we replace $z_w$ by any function $\psi\in V$. This is a consequence of the observation that the set $\{z_w : w\in V \}$ is dense in $V$. Indeed, given any $\psi\in V$ let $\{\eta_h\}_h\subset C^{\infty}_{\mathrm c}(0,L)$ be a sequence such that as $h\to\infty$
\begin{equation*}
\eta_h\to \psi'\quad\text{strongly in }L^2(0,L).
\end{equation*}
Setting $w_h=-\eta_h'$, one has that $w_h\in V$ and $z'_{w_h}=\eta_h$, hence $z_{w_h}\to\psi$ in $V$.
\end{proof}

\begin{lemma} \label{lem:mu''}
Let $u$ satisfy the assumptions of Proposition~\ref{prop:subdiff}.
The following identity holds:
\begin{equation} \label{nd366}
\int_0^L u'(x)w'(x)\,\de x + \int_0^L \partial_s W(x,u(x))w(x)\,\de x = \int_0^L z_{w}'(x)\mu'(x) \,\de x \qquad\text{for all } w\in V.
\end{equation}
\end{lemma}

\begin{proof}
By integration by parts we have, for all $w\in V$,
\begin{equation*}
\begin{split}
\int_0^L u'(x)&w'(x)\,\de x + \int_0^L \partial_s W(x,u(x))w(x)\,\de x \\
& = -\int_0^L u''(x)w(x)\,\de x+u'(L)w(L)-u'(0)w(0) +\int_0^L \partial_s W(x,u(x))w(x)\,\de x \\
& = \int_0^L \mu(x) w(x)\,\de x
\stackrel{\eqref{nd112ter}}{=} \scp{w}{\mu}_{V',V}
\stackrel{\eqref{nd108}}{=} \int_0^L z_{w}'(x)\mu'(x) \,\de x \,,
\end{split}
\end{equation*}
which is \eqref{nd366}.
\end{proof}

\section{Minimizing movements} \label{sect:MM}

Throughout this section, $q$ is a fixed element in $L^2(0,T;V)$, $\beta>0$, and $u_0\in V$ is a given initial datum.
We will exploit the variational structure of the equation
\begin{equation} \label{nd140}
\dot{u} = \bigl( -u'' + \partial_s W(x,u) \bigr)'' - \beta q',
\end{equation}
which was pointed out in the previous section, in order to construct, in Proposition~\ref{prop:compact}, a weak solution $u(x,t)$ to \eqref{nd140} satisfying the initial condition $u(x,0)=u_0(x)$, as the limit of a standard minimizing movement scheme. We name it \emph{variational solution} in Definition~\ref{def:varsol}, and in Proposition~\ref{prop:uni} we show that it is unique.

For $N\in\mathbb{N}$, let $\tau \coloneqq \frac{T}{N}$ be the time discretization step; we divide the time interval $[0,T]$ into subintervals $[(k-1)\tau,k\tau]$, $k=1,\ldots,N$, of length $\tau$.
We define recursively functions $u_\tau^k$ as follows: let $u_\tau^0(x) \coloneqq u_0(x)$ and, given $u_\tau^0,\ldots,u_\tau^{k-1}\in V$, $k\geq1$, let
\begin{equation}\label{nd141}
u_\tau^{k} \in \argmin \big\{
\Phi_\tau^k(v) : v\in V'
\big\},
\end{equation}
where, recalling the definition of $\Psi_q$ in \eqref{nd119},
\begin{equation*}
\Phi_\tau^k(v) \coloneqq \int_{k-1}^{k}\Psi_q(v,\tau t)\,\de t+\frac1{2\tau}\norma{u_\tau^{k-1}-v}_{V'}^2.
\end{equation*}
The functions $u_\tau^k$ are well-defined in view of the following lemma.

\begin{lemma} \label{lem:existMM}
The minimum problem \eqref{nd141} has a solution $u^k_\tau\in V$.
\end{lemma}

\begin{proof}
Let $\{v_n\}_n\subset V'$ be a minimizing sequence such that
\begin{equation}\label{nd144}
\Phi_\tau^k(v_n)\leq \inf_{v\in V'} \Phi_\tau^k(v)+\frac1n \leq \Phi_\tau^k(u_\tau^{k-1})+\frac1n =\int_{k-1}^{k} \Psi_q(u_\tau^{k-1},\tau t)\,\de t+\frac1n\,.
\end{equation}
By definition of $\Psi_q$, we have $v_n\in V$ for every $n$. In view of \eqref{nd122} we deduce the lower bound
\begin{equation}\label{nd145}
\begin{split}
\Phi_\tau^k(v_n)
\geq & \frac12\int_0^L|{v_n'(x)}|^2\de x-K_0L - \beta L^2\norma{v_n}_{L^2} \int_{k-1}^{k} \norma{q'(\tau t)}_{L^2} \,\de t
+\frac1{2\tau}\norma{u_\tau^{k-1}-v_n}_{V'}^2, \\
\end{split}
\end{equation}
where we have used \eqref{nd112bis} to estimate
\begin{equation*}
\big| (q'(\tau t),v_n)_{V'} \big| \leq \norma{q'(\tau t)}_{V'}\norma{v_n}_{V'} \leq L^2\norma{q'(\tau t)}_{L^2(0,L)}\norma{v_n}_{L^2(0,L)} .
\end{equation*}
Combining \eqref{nd144} and \eqref{nd145} and using Young's inequality we obtain for every $\eps>0$
\begin{equation}\label{nd148}
\begin{split}
\frac12\int_0^L \norm{v_n'(x)}^2\,\de x+\frac1{2\tau}\norma{u_\tau^{k-1}-v_n}_{V'}^2
& \leq K_0L + \beta L^2\eps\norma{v_n}_{L^2}^2 + \frac{\beta L^2}{4\eps} \int_{k-1}^{k} \norma{q'(\tau t)}_{L^2}^2\,\de t \\
& \quad +\int_{k-1}^{k} \Psi_q(u_\tau^{k-1},\tau t)\,\de t+\frac1n\,.
\end{split}
\end{equation}
Notice now that, since $q\in L^2(0,T;V)$, there holds
\begin{equation}\label{nd147}
\int_{k-1}^{k} \|q'(\tau t)\|_{L^2(0,L)}^2\,\de t
= \frac{1}{\tau}\int_{(k-1)\tau}^{k\tau}\|q'(t)\|_{L^2(0,L)}^2\,\de t\leq\frac{1}{\tau}\int_0^T \norma{q'(t)}^2_{L^2(0,L)}\,\de t.
\end{equation}
Moreover by choosing $\eps>0$ small enough, depending on $L$ and $\beta$, we have by Poincar\'e inequality (recalling that $v_n(0)=0$)
\begin{equation}\label{nd149}
\beta L^2\eps \norma{v_n}_{L^2(0,L)}^2
\leq\frac14\int_0^L \norm{v_n'(x)}^2\de x \,.
\end{equation}
In view of \eqref{nd147} and \eqref{nd149}, estimate \eqref{nd148} can be written as
\begin{equation*}\label{nd150}
\frac14\int_0^L \norm{v_n'(x)}^2\,\de x + \frac1{2\tau}\norma{u_\tau^{k-1}-v_n}_{V'}^2
\leq K_0L + \frac{\beta L^2}{4\eps\tau}\norma{q}_{L^2(0,T;V)}^2
+\int_{k-1}^{k} \Psi_q(u_\tau^{k-1},\tau t)\,\de t+\frac1n \,,
\end{equation*}
and yields the uniform bound $\sup_n \norma{v_n}_{V} <\infty$.
Therefore there exists $v\in V$ such that, up to subsequences, $v_n\wto v$ weakly in $V$ and uniformly in $[0,L]$; in particular, $v_n\to v$ in $V'$.
The minimality of $v$ in problem \eqref{nd141} follows then by the lower semicontinuity of $\Psi_q$ with respect to the above convergences.
\end{proof}

In the following lemma we compute the Euler-Lagrange equations satisfied by $u^k_{\tau}$ and we deduce its regularity properties.

\begin{lemma} \label{lem:ELEk}
Let $k\in\{1,\ldots,N\}$ and let $u_\tau^k$ be a minimizer of \eqref{nd141}. Then $u_\tau^k\in H^4(0,L)$ and
\begin{equation}\label{nd155}
-(u_\tau^k)'' + \displaystyle\partial_s W(\cdot,u_\tau^k(\cdot)) = -\beta\int_{k-1}^k z_{q'(\tau s)}\,\de s -z_{\bigl(\frac{u^k_{\tau}-u^{k-1}_{\tau}}{\tau}\bigr)},
\end{equation}
with boundary conditions
\begin{equation} \label{nd155bis}
\begin{cases}
u_\tau^k(0)=0, \\
\smallskip
(u_\tau^k)'(L)=0,\\
\smallskip
\bigl( -(u_{\tau}^k)''+\partial_s W(\cdot,u_\tau^k(\cdot)) \bigr)_{|_{x=0}}=0,\\
\smallskip
\bigl( -(u_{\tau}^k)''+\partial_s W(\cdot,u_\tau^k(\cdot)) \bigr)'_{|_{x=L}}=0.
\end{cases}
\end{equation}
In addition, 
\begin{align}\label{nd176}
\norma{u_\tau^k}_{H^3(0,L)}
& \leq C K \Bigl( 1 + \norma{u_\tau^k}_{H^1(0,L)} \Bigr) + \frac{C\beta}{\sqrt\tau}\norma{q}_{L^2((k-1)\tau,k\tau;V)} + \norma{\frac{u^k_{\tau}-u^{k-1}_{\tau}}{\tau}}_{V'}
\end{align}
and
\begin{equation}\label{nd177}
\begin{split}
\norma{(u_\tau^k)^{iv}}_{V'}
& \leq CK\biggl( K + K \|u_\tau^k\|_{H^1(0,L)}^2 + \frac{\beta}{\sqrt\tau}\norma{q}_{L^2((k-1)\tau,k\tau;V)} + \norma{\frac{u^k_{\tau}-u^{k-1}_{\tau}}{\tau}}_{V'} \biggr),
\end{split}
\end{equation}
where $C$ is a uniform constant depending only on $L$, and
\begin{equation} \label{nd180}
K=K\bigl(W,\|u^k_\tau\|_\infty\bigr) \coloneqq \| W \|_{C^3([0,L]\times S)}\,, \qquad S \coloneqq \bigl\{ s\in\R{} : |s|\leq\|u_\tau^k\|_\infty \bigr\}\,.
\end{equation}
\end{lemma}
\begin{proof}
We fix $w\in V$ and $\eps\in\R{}$.
By using $u_\tau^k+\eps w$ as a competitor in the minimum problem \eqref{nd141} solved by $u_\tau^k$ we have
\begin{equation*}\label{nd157}
\int_{k-1}^k \Psi_q(u_\tau^k,\tau s)\,\de s  +\frac1{2\tau}\norma{u_\tau^k-u_\tau^{k-1}}_{V'}^2
\leq \int_{k-1}^k \Psi_q(u_\tau^k+\eps w,\tau s)\,\de s +\frac1{2\tau}\norma{u_\tau^k+\eps w-u_\tau^{k-1}}_{V'}^2 ,
\end{equation*}
so that dividing by $\eps$ and letting $\eps\to 0$ we obtain that for every $w\in V$
\begin{equation}\label{nd159}
\int_0^L(u_\tau^k)'w'\,\de x+\int_0^L \partial_sW(x,u_\tau^k(x))w(x)\,\de x+\beta\int_{k-1}^k (q'(\tau s),w)_{V'}\,\de s + (\delta u_\tau^k,w)_{V'}=0,
\end{equation}
where we have defined $\delta u_\tau^k \coloneqq (u_\tau^k-u_\tau^{k-1})/\tau$ to simplify the notation.
We observe that by \eqref{nd110}, \eqref{nd108}, and \eqref{nd112ter}
\begin{equation}\label{nd160}
(\delta u_\tau^k,w)_{V'}
= \int_0^L z_{\delta u_\tau^k}'(x) z_{w}'(x)\,\de x
=\scp{w}{z_{\delta u_\tau^k}}_{V',V}
= \int_0^L z_{\delta u_\tau^k}(x)w(x)\,\de x \,,
\end{equation}
and similarly
\begin{equation}\label{nd161}
(q'(\tau s),w)_{V'}=\int_0^L z_{q'(\tau s)}(x)w(x)\,\de x \,.
\end{equation}
Notice that the map $t\mapsto z_{q'(t)}$ is measurable from $(0,T)$ to $V$, thanks to Lemma~\ref{lem:zq} applied to $q'$,
thus the map $t\mapsto\int_0^L z_{q'(t)}(x)w(x)\,\de x$ is measurable as well, and
\begin{equation}\label{nd169}
\int_0^T\bigg|\int_0^L z_{q'(t)}(x)w(x)\,\de x\bigg|\,\de t
=\int_0^T \big| (q'(t),w)_{V'} \big| \,\de t
\stackrel{\eqref{nd112bis}}{\leq} L^2 \norma{w}_{L^2(0,L)}\int_0^T \norma{q'(t)}_{L^2(0,L)}\de t<\infty.
\end{equation}
In view of estimate \eqref{nd169} and Fubini's theorem, the third term in \eqref{nd159} becomes
\begin{equation} \label{nd171}
\beta\int_{k-1}^k (q'(\tau s),w)_{V'}\,\de s
\stackrel{\eqref{nd161}}{=}
\beta\int_{k-1}^k \biggl( \int_0^L z_{q'(s\tau)}w\,\de x \biggr)\,\de s
= \beta\int_0^L \biggl(\int_{k-1}^k z_{q'(s\tau)}\,\de s\biggr)w\,\de x.
\end{equation}
Therefore, by combining \eqref{nd159}, \eqref{nd160} and \eqref{nd171}, we deduce that for every $w\in V$
\begin{equation*}
\begin{split}
\int_0^L(u_\tau^k)'w'\,\de x + \int_0^L &\partial_sW(x,u_\tau^k(x))w(x)\,\de x \\
& + \beta\int_0^L \biggl( \int_{k-1}^k z_{q'(s\tau)}(x)\,\de s \biggr) w(x)\,\de x
+ \int_0^L z_{\delta u_\tau^k}(x)w(x)\,\de x=0\,.
\end{split}
\end{equation*}
In particular, $u_\tau^k$ is a weak solution to the equation
\begin{equation} \label{nd163}
\begin{cases}
-(u_\tau^k)'' = - \partial_s W(\cdot,u_\tau^k(\cdot)) - \beta\int_{k-1}^k z_{q'(\tau s)}\,\de s -z_{\delta u_\tau^k}, \\
u_\tau^k(0)=0, \\
(u_\tau^k)'(L)=0.
\end{cases}
\end{equation}

We turn to the proof of the $H^4$-regularity of $u_\tau^k$, together with the bounds \eqref{nd176}--\eqref{nd177} and the last two conditions in \eqref{nd155bis}.
In the following estimates, $C$ will denote a positive constant depending possibly only on $L$, and which might change from line to line. Notice first that $u^k_\tau$ is a continuous function in $[0,L]$ and in particular $\|u_\tau^k\|_{L^\infty(0,L)}<\infty$. We can now estimate the right-hand side of \eqref{nd163} in $H^1(0,L)$ as follows. Using the regularity assumption on $W$ it is straightforward to obtain that
\begin{equation}\label{nd165}
\norma{\partial_s W(\cdot,u_\tau^k(\cdot))}_{H^1(0,L)} \leq C K \bigl( 1 + \|u_\tau^k\|_{H^1(0,L)} \bigr),
\end{equation}
where $K$ is the constant defined in \eqref{nd180}.
By \eqref{nd111}, $z_{q'(s)}\in H^2(0,L)$ with
\begin{equation*}
\norma{z_{q'(s)}}_{H^2(0,L)}\leq C\norma{q'(s)}_{L^2(0,L)} \leq C\norma{q(s)}_V,
\end{equation*}
and in turn
\begin{equation}\label{nd174}
\begin{split}
\bigg\| \int_{k-1}^k z_{q'(s\tau)}\,\de s\, \bigg\|_{H^2(0,L)}
& \leq \int_{k-1}^k \norma{z_{q'(s\tau)}}_{H^2(0,L)}\de s
\leq \biggl( \int_{k-1}^k \norma{z_{q'(s\tau)}}_{H^2(0,L)}^2\,\de s \biggr)^{\tfrac12} \\
& \leq C \biggl( \int_{k-1}^k \norma {q(s\tau)}_{V}^2\,\de s \biggr)^{\tfrac12}
\leq C \biggl( \frac{1}{\tau} \int_{(k-1)\tau}^{k\tau} \norma{q(s)}_{V}^2\,\de s \biggr)^{\tfrac12}.
\end{split}
\end{equation}
By \eqref{nd108}--\eqref{nd111} there holds
\begin{equation}\label{nd168}
\| z_{\delta u_\tau^k} \|_{H^1(0,L)} = \| \delta u_\tau^k \|_{V'},
\qquad
\| z_{\delta u_\tau^k}'' \|_{V'} = \| \delta u_\tau^k \|_{V'}.
\end{equation}
Combining \eqref{nd165}, \eqref{nd174}, and the first equality in \eqref{nd168}, we see that the right-hand side of the first equation in \eqref{nd163} is in $H^1(0,L)$, and therefore $u_\tau^k\in H^3(0,L)$ and the estimate \eqref{nd176} holds.

We can now further improve the regularity of $u^k_\tau$ by a bootstrap argument. Indeed, the regularity of $W$, \eqref{nd176}, and \eqref{nd165} yield
\begin{equation} \label{nd175}
\begin{split}
\norma{\partial_s W(\cdot,u_\tau^k(\cdot))}_{H^2(0,L)}
& \leq C K \Bigl( 1 + \|u_\tau^k\|_{H^1(0,L)}^2 + \|u_\tau^k\|_{H^2(0,L)}\Bigr) \\
& \leq CK\Bigl( K + K \|u_\tau^k\|_{H^1(0,L)}^2 + \frac{\beta}{\sqrt\tau}\norma{q}_{L^2((k-1)\tau,k\tau;V)} + \norma{\delta u_\tau^k}_{V'} \Bigr).
\end{split}
\end{equation}
Therefore, by \eqref{nd175}, \eqref{nd174}, and \eqref{nd111}, the right-hand side of the first equation in \eqref{nd163} belongs to $H^2(0,L)$, so that $u_\tau^k\in H^4(0,L)$. Furthermore, we can estimate the $V'$-norm of the fourth derivative of $u_\tau^k$ by the $V'$-norm of the second derivative of the right-hand side of \eqref{nd163}: combining \eqref{nd175}, \eqref{nd174}, and the second equation in \eqref{nd168}, we obtain \eqref{nd177}.

Finally, it is straightforward to check that
$$
-\Big(\int_{k-1}^k z_{q'(t\tau)}\,\de t\Big)''=\int_{k-1}^k q'(t\tau)\,\de t \,,
$$
and since the function $\int_{k-1}^k z_{q'(t\tau)}\,\de t$ vanishes in $x=0$ and has first derivative vanishing in $x=L$, we conclude by \eqref{nd111} that
\begin{equation}\label{nd177bis}
\int_{k-1}^k z_{q'(t\tau)}\,\de t=z_{\int_{k-1}^kq'(t\tau)\,\de t}.
\end{equation}
Thus, the regularity of $u_{\tau}^k$ and \eqref{nd155} yield the last two conditions in \eqref{nd155bis}.
\end{proof}

We next introduce the piecewise affine and the piecewise constant interpolants of the functions $u_\tau^k$:
\begin{subequations}\label{nd153}
\begin{eqnarray}
u_\tau(t) & \!\!\!\! \coloneqq & \!\!\!\! u_\tau^{k-1}+\Big(\frac{t-(k-1)\tau}\tau\Big)(u_\tau^k-u_\tau^{k-1}),\quad\text{for $(k-1)\tau\leq t\leq k\tau$,} \label{nd153a} \\
\tilde{u}_\tau(t) & \!\!\!\! \coloneqq & \!\!\!\! u_\tau^{k},\quad\text{for } (k-1)\tau < t \leq k\tau. \label{nd153b}
\end{eqnarray}
\end{subequations}
Notice that $u_\tau(0,t)=0$ for all $t$, and, in particular, $u_\tau(t)\in V$ for all $t$. In addition, $u_\tau(x,0)=\tilde{u}_\tau(x,0)=u_0(x)$, and 
\begin{equation*}
\dot u_\tau(x,t)=\frac{u_\tau^k-u_\tau^{k-1}}\tau\qquad\text{for $(k-1)\tau<t<k\tau$.}
\end{equation*}
By Lemma \ref{lem:ELEk} we deduce the following uniform estimates.

\begin{lemma} \label{lem:unifest}
Let $u_\tau$ and $\tilde{u}_\tau$ be defined as in \eqref{nd153a} and \eqref{nd153b}, respectively.
There exist constants $M_1$ and $M_2$, dependent on $\norma{q}_{L^2(0,T;V)}$, $\beta$, $\mathcal{E}(u_0)$, $L$, $T$ and on the potential $W$, such that
\begin{align}
\sup_\tau \norma{u_\tau}_{C^0([0,T];V)} \leq{} & M_1, \label{nd189} \\
\sup_\tau \norma{\dot u_\tau}_{L^2(0,T;V')} \leq{} & M_1, \label{nd190} \\
\sup_\tau  \norma{\tilde{u}_\tau}_{L^2(0,T;H^3(0,L))} \leq{} & M_2, \label{nd191} \\
\sup_\tau  \| \tilde{u}_\tau^{(iv)} \|_{L^2(0,T;V')} \leq{} & M_2. \label{nd191ter}
\end{align}

\end{lemma}

\begin{proof}
By minimality of $u_\tau^k$ and using $u_\tau^{k-1}$ as a competitor in the minimum problem \eqref{nd141}, we have $\Phi_\tau^k(u_\tau^k) \leq  \Phi_\tau^k(u_\tau^{k-1})$, or explicitly
\begin{equation*}
\begin{split}
\frac12 \int_0^L & \norm{(u_\tau^k)'(x)}^2\,\de x + \int_0^L W(x,u_\tau^k(x))\,\de x+\beta\int_{k-1}^k (q'(t\tau),u_\tau^k)_{V'}\,\de t+\frac12\int_{(k-1)\tau}^{k\tau} \norma{\dot u_\tau(t)}_{V'}^2\de t \\
& \leq \frac12\int_0^L \norm{(u_\tau^{k-1})'(x)}^2\,\de x + \int_0^L W(x,u_\tau^{k-1}(x))\,\de x+\beta\int_{k-1}^k (q'(t\tau),u_\tau^{k-1})_{V'}\,\de t.
\end{split}
\end{equation*}
Iterating this estimate it follows that
\begin{equation*}
\begin{split}
& \frac12\int_0^L \norm{(u_\tau^k)'(x)}^2\,\de x + \int_0^L W(x,u_\tau^{k}(x))\,\de x + \beta\sum_{j=1}^{k}\int_{j-1}^j (q'(t\tau),u_\tau^j-u_\tau^{j-1})_{V'}\,\de t \\
& \qquad\qquad\qquad + \frac12\int_{0}^{k\tau} \norma{\dot u_\tau(t)}_{V'}^2\de t\leq \frac12\int_0^L \norm{u_0'(x)}^2\de x + \int_0^L W(x,u_0(x))\,\de x=\mathcal{E}(u_0),
\end{split}
\end{equation*}
namely, by \eqref{nd122},
\begin{equation}\label{nd183}
\frac12\int_0^L \norm{(u_\tau^k)'(x)}^2\,\de x
+\int_{0}^{k\tau} \bigl( \beta q'(t)+\textstyle\frac12\dot u_\tau(t),\dot u_\tau(t) \bigr)_{V'}\,\de t \leq \mathcal{E}(u_0) + K_0L,
\end{equation}
for every $k=1,\ldots,N$. Young's inequality implies
\begin{equation*}
\begin{split}
\int_0^{k\tau} & \bigl( \beta q'(t) +{\textstyle\frac12}\dot u_\tau(t),\dot u_\tau(t) \bigr)_{V'}\,\de t
\geq -\beta\int_0^{k\tau} \norma{q'(t)}_{V'}\norma{\dot u_\tau(t)}_{V'}\,\de t + \frac12\int_0^{k\tau} \norma{\dot u_\tau(t)}_{V'}^2\,\de t \\
& \geq -\frac14\int_0^{k\tau} \norma{\dot u_\tau(t)}_{V'}^2\de t - \beta^2\int_0^{k\tau} \norma{q'(t)}_{V'}^2\de t + \frac12\int_0^{k\tau} \norma{\dot u_\tau(t)}_{V'}^2\,\de t \\
& = \frac14\int_0^{k\tau} \norma{\dot u_\tau(t)}_{V'}^2\de t-\beta^2\int_0^{k\tau} \norma{q'(t)}_{V'}^2\de t,
\end{split}
\end{equation*}
which, combined with \eqref{nd183}, yields, in view of \eqref{nd112bis},
\begin{equation*}
\frac12\int_0^L \norm{(u_\tau^k)'(x)}^2\,\de x+\frac14\int_{0}^{k\tau} \norma{\dot u_\tau(t)}_{V'}^2\de t\leq L^2\beta^2\int_0^T \norma{q'(t)}_{L^2(0,L)}^2\de t+\mathcal{E}(u_0) + K_0L,
\end{equation*}
where the integral on the right-hand side is finite since $q\in L^2(0,T;V)$.
It follows that
\begin{equation}\label{nd186}
\sup_{k,\tau}\norma{u_\tau^k}_{H^1(0,L)} \leq M_1, \qquad \sup_\tau\norma{\dot u_\tau}_{L^2(0,T;V')} \leq M_1,
\end{equation}
with $M_1= 2L\beta \norma{q}_{L^2(0,T;V)} + 2\sqrt{\mathcal{E}(u_0) + K_0L}$. Therefore estimates \eqref{nd189} and \eqref{nd190} are a direct consequence of \eqref{nd186}, by definition of $u_\tau$.

Notice now that the constant $K$ defined in \eqref{nd180}, which in principle depends on the $L^\infty$-norm of $u_\tau^k$, is uniformly bounded with respect to $\tau$ and $k$ by a constant $M_2$ as in the statement, in view of \eqref{nd186}. Then property \eqref{nd191} follows by integrating \eqref{nd176} in the interval $((k-1)\tau,k\tau)$, summing over $k$, and using \eqref{nd189}--\eqref{nd190}. Similarly, \eqref{nd191ter} follows by integrating \eqref{nd177} in $((k-1)\tau,k\tau)$ and summing over $k$.
\end{proof}

The uniform estimates obtained in Lemma~\ref{lem:unifest} allow us to extract a subsequence of $\{u_\tau\}$ converging, as $\tau\to0$, to a limit function $u_\infty$. In turn, by passing to the limit in the equation \eqref{nd155} satisfied by $u_\tau^k$, we deduce that $u_\infty$ actually solves problem \eqref{nd140} in a weak sense. This is the content of the following proposition.

\begin{proposition} \label{prop:compact}
There exists a map $u_\infty\in L^2(0,T;H^3(0,L))\cap H^1(0,T;V')$ 
such that, up to the extraction of a (not relabeled) subsequence, the following convergences hold true as $\tau\to0$:
\begin{enumerate}
	\item\label{itemcomp1} $u_\tau\to u_\infty$ in $C^0([0,T];V')$,
	\item\label{itemcomp2} $u_\tau(t)\wto u_\infty(t)$ weakly in $H^1(0,L)$ for every $t\in[0,T]$,
	\item\label{itemcomp3} $u_\tau \stackrel{*}{\wto} u_\infty$ weakly* in $L^\infty(0,T;V)$,
	\item\label{itemcomp4} $\tilde{u}_\tau \wto u_\infty$ weakly in $L^2(0,T;H^3(0,L))$.
\end{enumerate}
In addition, $u_\infty$ solves for almost every $(x,t)\in(0,L)\times(0,T)$
\begin{equation} \label{nd201}
-u_\infty''(x,t) + \partial_sW(x,u_\infty(x,t)) = - z_{\dot{u}_\infty(t)}(x) - \beta z_{q'(t)}(x)
\end{equation}
with $u_\infty(x,0)=u_0(x)$, $u_\infty(0,t)=u_\infty'(L,t)=0$.
\end{proposition}

\begin{proof}
We divide the proof of the proposition into two steps.

\smallskip
\noindent
\textit{Step 1: compactness.}
In view of \eqref{nd190} we obtain for all $0\leq s \leq t \leq T$
\begin{equation}\label{nd192}
\begin{split}
\sup_\tau \norma{u_\tau(t)-u_\tau(s)}_{V'}^2\leq{} & \sup_\tau \left(\int_s^t \norma{\dot u_\tau(r)}_{V'}\,\de r\right)^2 \\
\leq{} & (t-s)\sup_\tau \int_s^t\norma{\dot u_\tau(r)}_{V'}^2\de r\leq M_{1}^2(t-s),
\end{split}
\end{equation}
which shows that the family of functions $\{u_\tau\}_\tau$ is uniformly equicontinuous from $[0,T]$ to $V'$.
In addition, $\{u_\tau(t)\}_\tau$ is relatively compact in $V'$ for all $t\in[0,T]$: indeed, by \eqref{nd189} the functions $u_\tau(t)$ are equibounded in $V$, and thus up to subsequences they converge strongly in $L^2(0,L)$ and, in turn, in $V'$.
Hence by Ascoli-Arzel\`a Theorem we can extract a (not relabeled) subsequence such that \ref{itemcomp1} holds.
Since the sequence $(u_\tau)_\tau$ is uniformly bounded in $H^1(0,T;V')$ by \eqref{nd190}, it follows that 
\begin{equation} \label{nd198}
\dot{u}_\tau \wto \dot{u}_\infty \qquad\text{weakly in }L^2(0,T;V').
\end{equation}
and $u_\infty\in H^1(0,T;V')$. 
Properties~\ref{itemcomp2} and \ref{itemcomp3} are consequences of the uniform estimate \eqref{nd189}. 

Notice that for all $t\in(0,T]$ we have $\tilde{u}_\tau(t)=u_\tau^{k}=u_\tau(k\tau)$ for some $k\in\{1,\ldots,N\}$, hence by \eqref{nd192} and \ref{itemcomp1}
\begin{equation*}
\begin{split}
\norma{\tilde{u}_\tau(t)-u_\infty(t)}_{V'}
& \leq\norma{u_\tau(k\tau)-u_\tau(t)}_{V'}+\norma{u_\tau(t)-u_\infty(t)}_{V'} \\
& \leq M_1\sqrt\tau+\norma{u_\tau(t)-u_\infty(t)}_{V'}\to0,
\end{split}
\end{equation*}
that is,
\begin{equation*}
\tilde{u}_\tau(t)\to u_\infty(t) \qquad\text{strongly in $V'$, for every }t\in[0,T].
\end{equation*}
This identifies the pointwise limit of $\tilde{u}_\tau$; then  estimate \eqref{nd191} implies \ref{itemcomp4} and yields the regularity $u_\infty\in L^2(0,T;H^3(0,L))$.

\smallskip
\noindent
\textit{Step 2: derivation of the equation for $u_\infty$.} We turn to the proof of \eqref{nd201}. We can rewrite the discrete equation \eqref{nd155} using the interpolants $u_\tau$, $\tilde{u}_\tau$ in the form
\begin{equation} \label{nd199bis}
\begin{split}
-\tilde{u}_\tau''(x,t) + \partial_sW(x,&\tilde{u}_\tau(x,t)) = -z_{\dot{u}_\tau(t)}(x) - \beta \sum_{k=1}^N\biggl(\int_{k-1}^k z_{q'(\tau s)}\de s\biggr)\chi_{((k-1)\tau,k\tau)}(t)
\end{split}
\end{equation}
for almost every $(x,t)\in(0,L)\times(0,T)$. Notice that by \eqref{nd177bis} the last term in \eqref{nd199bis} can be written as $\beta z_{r_\tau(t)}$, where
\begin{equation*}
r_\tau(t) \coloneqq \frac{1}{\tau}\sum_{k=1}^N \chi_{((k-1)\tau,k\tau)}(t) \int_{(k-1)\tau}^{k\tau} q'(s)\,\de s\,.
\end{equation*}
Then by testing \eqref{nd199bis} with an arbitrary function $\varphi\in C^\infty([0,L]\times[0,T])$ we have
\begin{multline} \label{nd199}
\int_0^T\int_0^L \bigl( -\tilde{u}_\tau''(x,t) + \partial_sW(x,\tilde{u}_\tau(x,t)) \bigr) \varphi(x,t) \,\de x\de t \\
= - \int_0^T\int_0^L z_{\dot{u}_\tau(t)}(x)\varphi(x,t)\, \de x \de t
- \beta \int_0^T \int_0^L z_{r_\tau(t)}(x)\varphi(x,t)\,\de x \de t \,.
\end{multline}

In order to pass to the limit as $\tau\to0$ in \eqref{nd199}, observe first that by \ref{itemcomp4} we have $\tilde{u}_\tau''\wto u_\infty''$ weakly in $L^2(0,T,H^1(0,L))$. By \eqref{nd186}, for every fixed $t\in[0,T]$ the sequence $(\tilde{u}_\tau(t))_\tau$ is uniformly bounded in $V$, so that we have also $\tilde{u}_\tau(t)\to u_\infty(t)$ uniformly in $[0,L]$ and, in turn, $\partial_sW(x,\tilde{u}_\tau(x,t))\to\partial_sW(x,u_\infty(x,t))$ for every $(x,t)$. By \eqref{nd198} there holds $z_{\dot{u}_\tau}\wto z_{\dot{u}_\infty}$ weakly in $L^2(0,T;V)$. Finally, Lemma~\ref{lem:q} yields that $r_\tau\wto q'$ weakly in $L^2(0,T;L^2(0,L))$, so that $z_{r_\tau(t)}\wto z_{q'(t)}$ weakly in $L^2(0,T;V)$.

In view of all these convergences we can then pass to the limit as $\tau\to0$ in \eqref{nd199}:
\begin{multline*}
\int_0^T\int_0^L \bigl( -u_\infty''(x,t) + \partial_sW(x,u_\infty(x,t)) \bigr) \varphi(x,t) \,\de x\de t \\
= - \int_0^T\int_0^L z_{\dot{u}_\infty(t)}(x)\varphi(x,t)\, \de x \de t
- \beta \int_0^T \int_0^L z_{q'(t)}(x)\varphi(x,t)\,\de x \de t \,,
\end{multline*}
and as the test function $\varphi\in C^\infty([0,L]\times[0,T])$ is arbitrary we conclude that the equation \eqref{nd201} holds pointwise almost everywhere in $(0,L)\times(0,T)$.
\end{proof}

\begin{defin} \label{def:varsol}
We say that $u$ is a \emph{variational solution} to \eqref{nd140} with initial datum $u_0\in V$ if $u$ is the limit of a subsequence of minimizing movements as in Proposition~\ref{prop:compact}.
\end{defin}

The following proposition implies in particular the uniqueness of the variational solution for every given $q\in L^2(0,T;V)$, since, by Proposition~\ref{prop:compact}, every variational solution satisfies the regularity assumptions in Proposition \ref{prop:uni}.

\begin{proposition} \label{prop:uni}
Let $q\in L^2(0,T;V)$ and let $u_0\in V$. Then there exists a unique function $u_q\in L^2(0,T;H^3(0,L))\cap H^1(0,T;V')$ such that for almost every $(x,t)\in(0,L)\times(0,T)$
\begin{equation} \label{nd200}
-u_q''(x,t) + \partial_sW(x,u_q(x,t)) = - z_{\dot{u}_q(t) + \beta q'(t)}(x)
\end{equation}
with $u_q(x,0)=u_0(x)$ and $u_q(0,t)=u_q'(L,t)=0$.
\end{proposition}

\begin{proof}
Let $u_1,u_2$ be two solutions of \eqref{nd200}, corresponding to the same $q$, satisfying the assumptions in the statement, and let $w\coloneqq u_1-u_2$. Then $w$ solves the equation
\begin{equation}\label{nd202}
-w''(x,t) + \partial_sW(x,u_1(x,t)) - \partial_sW(x,u_2(x,t)) = -z_{\dot{w}(t)}(x)
\end{equation}
with $w(x,0)=0$, $w(0,t)=w'(L,t)=0$.
Multiplying \eqref{nd202} by $z_{w(t)}$ and integrating in $(0,L)$ we obtain that for a.e.\ $t\in(0,T)$
\begin{equation} \label{nd203}
\begin{split}
\frac12&\frac{\de}{\de t}\norma{z_{w(t)}}^2_{L^2(0,L)}
= \int_0^L z_{\dot{w}(t)}(x)z_{w(t)}(x)\,\de x \\
&= -\int_0^L \bigl[ -w''(x,t) + \partial_sW(x,u_1(x,t)) - \partial_sW(x,u_2(x,t)) \bigr] z_{w(t)}(x)\,\de x \\
& \leq \int_0^L w(x,t)z_{w(t)}''(x)\,\de x + \norma{z_{w(t)}}_{L^2(0,L)} \norma{\partial_sW(\cdot,u_1(\cdot,t)) - \partial_sW(\cdot,u_2(\cdot,t))}_{L^2(0,L)},
\end{split}
\end{equation}
where the first equality is a consequence of Lemma~\ref{lem:zpunto} and of \cite[Theorem~8.60]{Leo}, and we integrated by parts twice in the last passage. Observe now that, using the smoothness \eqref{nd121} of the potential $W$ and the fact that $u_i\in L^\infty(0,T;L^\infty(0,L))$, $i=1,2$ (see Remark~\ref{rmk:interpolation}), we can find a constant $C>0$, independent of $t\in[0,T]$, such that
\begin{equation} \label{nd204}
\norma{\partial_sW(\cdot,u_1(\cdot,t)) - \partial_sW(\cdot,u_2(\cdot,t))}_{L^2(0,L)}
\leq C\norma{w(t)}_{L^2(0,L)}.
\end{equation}
Inserting \eqref{nd204} in \eqref{nd203}, and recalling \eqref{nd111}, we obtain
\begin{equation*}
\begin{split}
\frac{\de}{\de t}\norma{z_{w(t)}}^2_{L^2(0,L)}
& \leq - 2\norma{w(t)}_{L^2(0,L)}^2 + 2C \norma{z_{w(t)}}_{L^2(0,L)}\norma{w(t)}_{L^2(0,L)} \\
& \leq - \norma{w(t)}_{L^2(0,L)}^2 + C^2\norma{z_{w(t)}}^2_{L^2(0,L)},
\end{split}
\end{equation*}
where we used Young's inequality in the second passage.
Since $w(t)$ (and therefore $z_{w(t)}$) vanishes for $t=0$, by Gr\"onwall's Lemma we conclude that $z_{w(t)}=0$ and, in turn, $w(t)=0$ for all $t\in [0,T]$.
\end{proof}

\section{Proof of Theorem~\ref{thm:existence}} \label{sect:fixedpoint}

In this section we give the main argument for the proof of Theorem~\ref{thm:existence}. We divide the proof into three parts: in the first, we show the existence of a weak solution by means of a fixed point argument, and the uniqueness of such solution; subsequently we prove the energy equality \eqref{nd360}, and finally we show the continuous dependence on initial data and the estimate \eqref{S1}.

\subsection{Fixed point argument}
Let $u_0\in V$ be a given initial datum. Thanks to the results in the previous section, we can define a map $S\colon L^2(0,T;V)\to L^2(0,T;V)$ as $S(q)\coloneqq u_q$ for every $q\in L^2(0,T;V)$, where $u_q$ is the unique variational solution to \eqref{nd140} with initial datum $u_0$ constructed in Proposition~\ref{prop:compact}. Notice that the uniqueness of the variational solution is a consequence of Proposition~\ref{prop:uni}, hence the map $S$ is well defined. We now prove the existence of a weak solution to \eqref{sub2} by showing that the map $S$ has a fixed point.

\smallskip
\noindent
\textit{Step 1: continuity and compactness of $S$.}
Let $q_n\to q$ in $L^2(0,T;V)$. The estimates proved in Lemma~\ref{lem:unifest} on $u_\tau^{(n)}$ and on  $\tilde{u}_\tau^{(n)}$ (where the superscript $(n)$ indicates that these functions are the piecewise affine and piecewise constant interpolants corresponding to $q_n$, respectively) are uniform in $n$, since they depend on $n$ only through the $L^2(0,T;V)$ norm of $q_n$. Hence, passing to the limit as $\tau\to 0$ as in Proposition~\ref{prop:compact}, we obtain that the sequence $(S(q_n))_n$ is uniformly bounded in $L^2(0,T;H^3(0,L))\cap H^1(0,T;V')$, and converges weakly to $S(q)$ in the same spaces.

By Aubin-Lions-Simon Theorem (see \cite[Theorem~8.62]{Leo}), $L^2(0,T;H^3(0,L))\cap H^1(0,T;V')$ is compactly imbedded in $L^2(0,T;H^2(0,L))$, so that $S(q_n)\to S(q)$ strongly in $L^2(0,T;V)$. This proves the continuity of the map $S$ in $L^2(0,T;V)$.

Similarly, if $(q_n)_n$ is a bounded sequence in $L^2(0,T;V)$, we can repeat the previous argument to deduce that $(S(q_n))_n$ is uniformly bounded in $L^2(0,T;H^3(0,L))\cap H^1(0,T;V')$. In turn, as before the compact embedding given by Aubin Theorem yields that a subsequence of $S(q_n)$ converges strongly in $L^2(0,T;V)$, thus proving the compactness of the map $S$.

\smallskip
\noindent
\textit{Step 2: fixed point.}
In order to apply Shaefer's Fixed Point Theorem to the map $S$, we still need to check that the set
\begin{equation*}
\Lambda \coloneqq \Bigl\{ q\in L^2(0,T;V) \,:\, q=\lambda S(q) \text{ for some } 0\leq \lambda \leq 1 \Bigr\}
\end{equation*}
is bounded in $L^2(0,T;V)$. Let $q\in\Lambda$; by using the identity $u_q = q/\lambda$ 
in equation \eqref{nd201} solved by $u_q$ we obtain
\begin{align} \label{nd301}
-q''(t) + \lambda\,\partial_sW(\cdot,\textstyle\frac{1}{\lambda}q(\cdot,t)) = -z_{\dot{q}(t) + \beta\lambda q'(t)}\,.
\end{align}
Multiplying \eqref{nd301} by $q$, and integrating in $(0,L)$, we have
\begin{align} \label{nd302}
\int_0^L |q'(x,t)|^2\,\de x + \lambda \int_0^L \partial_sW(x,{\textstyle\frac{1}{\lambda}}q(x,t))q(x,t)\,\de x = - \int_0^L z_{\dot{q}(t) + \beta\lambda q'(t)}(x) q(x,t)\,\de x\,.
\end{align}
Notice that we integrated by parts in the first term, with the boundary terms vanishing as $q(0,t)=0$, $q'(L,T)=\lambda u'_q(L,t)=0$. Using the equality $-z_{q(t)}''=q(t)$, we can integrate by parts also the right-hand side of \eqref{nd302}, and we obtain
\begin{align*}
&- \int_0^L z_{\dot{q}(t) + \beta\lambda q'(t)}(x) q(x,t)\,\de x
=\int_0^L z_{\dot{q}(t) + \beta\lambda q'(t)}(x) z''_{q(t)}(x)\,\de x\,\\
&\quad=-\int_0^L z'_{\dot{q}(t) + \beta\lambda q'(t)}(x) z'_{q(t)}(x)\,\de x
=-\int_0^L z'_{\dot{q}(t)}(x) z'_{q(t)}(x)\,\de x-\beta\lambda\int_0^L q'(x,t)z_{q(t)}(x)\,\de x\,,
\end{align*}
which inserted in \eqref{nd302} yields the equality
\begin{equation}\label{nd303}
\begin{split}
\int_0^L z_{\dot{q}(t)}'(x)z_{q(t)}'(x)&\,\de x + \int_0^L |q'(x,t)|^2\,\de x \\
& = - \lambda \int_0^L \partial_sW(x,{\textstyle\frac{1}{\lambda}}q(x,t))q(x,t)\,\de x
- \beta\lambda \int_0^L q'(x,t)z_{q(t)}(x)\,\de x\,.
\end{split}
\end{equation}
Notice that, since $q=\lambda u_q$, we have in particular $q\in H^1(0,T;V')$. Therefore by Lemma~\ref{lem:zpunto} the first integral on the left-hand side of \eqref{nd303} can be written as $\frac12\frac{\de}{\de t}\int_0^L |z_{q(t)}'(x)|^2\,\de x$. Moreover, recalling \eqref{nd122} we have the estimate
\begin{equation} \label{nd304}
- \lambda \int_0^L \partial_sW(x,{\textstyle\frac{1}{\lambda}}q(x,t))q(x,t)\,\de x
\leq \lambda^2(K_0+K_1)L \leq (K_0+K_1)L \,,
\end{equation}
while by Young and Poincar\'e inequalities
\begin{equation} \label{nd304bis}
- \beta\lambda \int_0^L q'(x,t)z_{q(t)}(x)\,\de x
\leq \frac12\norma{q(t)}_V^2 + \frac{\beta^2L^2}{2}\norma{z_{q(t)}}_{V}^2 .
\end{equation}
Hence, inserting \eqref{nd304} and \eqref{nd304bis} in \eqref{nd303} we obtain
\begin{align} \label{nd305}
\frac{\de}{\de t} \norma{z_{q(t)}}_V^2 + \norma{q(t)}_V^2
\leq 2(K_0+K_1)L + \beta^2L^2\norma{z_{q(t)}}_{V}^2 .
\end{align}
Gr\"{o}nwall's Lemma then yields
\begin{align} \label{nd306}
\norma{z_{q(t)}}_V^2 \leq  \bigg( \norma{z_{q(0)}}_V^2 + \frac{2(K_0+K_1)}{\beta^2 L} \bigg) e^{\beta^2L^2T} \leq\overline{C}
\qquad \text{for every }t\in[0,T],
\end{align}
where the constant $\overline{C}$ on the right-hand side is independent of $q\in\Lambda$, as $q(0)=\lambda u_q(0)=\lambda u_0$. 
Integrating \eqref{nd305} in $[0,T]$ and using \eqref{nd306} we deduce the estimate
\begin{equation*}
\norma{z_{q(T)}}_V^2 - \norma{z_{q(0)}}_V^2 + \norma{q}_{L^2(0,T;V)}^2 \leq 2 (K_0+K_1)LT + \beta^2L^2\overline{C}T \,,
\end{equation*}
from which we finally get that $\norma{q}_{L^2(0,T;V)}$ is uniformly bounded for $q\in \Lambda$. Hence the set $\Lambda$ is bounded, as claimed.

We can then invoke Shaefer's Fixed Point Theorem to obtain the existence of an element $q\in L^2(0,T;V)$ such that $q=S(q)=u_q$.
By construction, $u_q$ belongs to the space $L^2(0,T;H^3(0,L))\cap H^1(0,T;V')$, solves equation \eqref{nd201}, that is
\begin{equation*}
-u_q''(x,t) + \partial_sW(x,u_q(x,t)) = - z_{\dot{u}_q(t)}(x) - \beta z_{u_q'(t)}(x)
\end{equation*}
for almost every $(x,t)\in(0,L)\times(0,T)$, and attains the boundary and initial conditions $u_q(0,t)=u_q'(L,t)=0$, $u_q(x,0)=u_0(x)$. The associated function $\mu=-u_q''+\partial_sW(x,u_q)$ satisfies the identity
\begin{equation*}
\mu(t) = - z_{\dot{u}_q(t)} - \beta z_{u_q'(t)},
\end{equation*}
and since by Lemma~\ref{lem:zq} the maps $t\mapsto z_{\dot{u}_q(t)}$, $t\mapsto z_{u_q'(t)}$ belong to $L^2(0,T;V)$, we also have $\mu\in L^2(0,T;V)$. Therefore, in view of Proposition~\ref{prop:weaksol}, all the conditions in Definition~\ref{def:weaksol} are satisfied and $u_q$ is a weak solution to the problem \eqref{sub2} in $[0,T]$ corresponding to the initial datum $u_0$.

\smallskip
\noindent
\textit{Step 3: uniqueness.}
Let $u_1,u_2$ be two weak solutions to \eqref{sub2}, corresponding to the same initial datum, satisfying the assumptions in the statement, and let $w\coloneqq u_1-u_2$. Then $w$ solves the equation
\begin{equation}\label{nd308}
-w''(x,t) + \partial_sW(x,u_1(x,t)) - \partial_sW(x,u_2(x,t)) = -z_{\dot{w}(t)}(x) - \beta z_{w'(t)}(x)
\end{equation}
with $w(x,0)=0$, $w(0,t)=w'(L,t)=0$.
Multiplying \eqref{nd308} by $w(t)$ and integrating the equation in $(0,L)$ we obtain
\begin{equation} \label{nd309}
\begin{split}
\int_0^L &|w'(x,t)|^2\,\de x + \int_0^L \bigl[ \partial_sW(x,u_1(x,t)) - \partial_sW(x,u_2(x,t)) \bigr] w(x,t)\,\de x \\
&= -\int_0^L \bigl( z_{\dot{w}(t)}(x) +\beta z_{w'(t)}(x) \bigr) w(x,t)\,\de x
= \int_0^L \bigl( z_{\dot{w}(t)}(x) + \beta z_{w'(t)}(x) \bigr) z''_{w(t)}(x)\,\de x \\
& = -\int_0^L z'_{\dot{w}(t)}(x)z'_{w(t)}(x)\,\de x -\beta\int_0^L z_{w(t)}(x)w'(x,t)\,\de x \\
& = -\frac12\frac{\de}{\de t}\norma{z_{w(t)}}_V^2 -\beta\int_0^L z_{w(t)}(x)w'(x,t)\,\de x\,,
\end{split}
\end{equation}
where we repeatedly integrated by parts and we used Lemma~\ref{lem:zpunto} in the last passage, as $w\in H^1(0,T;V')$. By the regularity \eqref{nd121} of the potential $W$, and observing that
\begin{equation*}
\sup_{t\in(0,T)} \|u_i(t)\|_{L^\infty(0,L)} <\infty, \qquad i=1,2
\end{equation*}
by the condition $u_i\in L^\infty(0,T;V)$, we can find a constant $C>0$, independent of $t\in[0,T]$ and $x\in[0,L]$, such that
\begin{equation*}
\big| \partial_sW(x,u_1(x,t)) - \partial_sW(x,u_2(x,t)) \big|
\leq C |w(x,t)| \,.
\end{equation*}
This estimates gives in turn
\begin{equation} \label{nd310}
\begin{split}
\bigg| \int_0^L \bigl[ \partial_sW(x,u_1(x,t)) - \partial_sW(x,u_2(x,t)) \bigr] w(x,t)\,\de x \bigg|
& \leq C\int_0^L |w(x,t)|^2\de x \\
& \leq C \norma{w(t)}_{V'}\norma{w(t)}_V.
\end{split}
\end{equation}
Hence, by combining  \eqref{nd309} and \eqref{nd310} we have
\begin{equation*}
\begin{split}
\frac12\frac{\de}{\de t}\norma{z_{w(t)}}_V^2
& \leq -\norma{w(t)}_V^2 + \beta\norma{z_{w(t)}}_{L^2(0,L)}\norma{w(t)}_V + C\norma{w(t)}_{V'}\norma{w(t)}_V \\
& \stackrel{\eqref{nd108}}{\leq} -\norma{w(t)}_V^2 + C\norma{z_{w(t)}}_V\norma{w(t)}_V \\
& \leq -\frac12\norma{w(t)}_V^2 + C\norma{z_{w(t)}}_V^2\,,
\end{split}
\end{equation*}
where we used Young's inequality in the last line.
Since $w(t)$ vanishes for $t=0$, by Gr\"onwall's Lemma we conclude that $z_{w(t)}\equiv0$ and, in turn, that $w(t)\equiv0$ for all $t$. This concludes the proof of uniqueness of the weak solution.

\subsection{Energy equality}
We now turn to the proof of the energy equality \eqref{nd360}. 
By taking $\psi=\mu(t)$ as test function in the weak formulation \eqref{nd358} (which is admissible since $\mu(t)\in V$ for almost every $t\in(0,T)$), we find
\begin{equation*}
\scp{\dot u(t)}{\mu(t)}_{V',V} = - \| \mu(t) \|^2_V -\beta\int_0^L u'(x,t)\mu(x,t)\,\de x \,.
\end{equation*}
To conclude, we only need to show that for almost every $t\in(0,T)$
\begin{equation}\label{nd362}
\frac{\de}{\de t}\mathcal{E}(u(t))= \scp{\dot u(t)}{\mu(t)}_{V',V} \,.
\end{equation}
Given $h>0$ and $\lambda\in(0,1)$, defining $u_h^\lambda(x,t)\coloneqq \lambda u(x,t)+(1-\lambda)u(x,t-h)$, we start by considering the increment of the energy from $t-h$ to $t$:
\begin{equation*}
\begin{split}
\mathcal{E}(u(t))-\mathcal{E}(u(t-h))
={} & \frac12\int_0^L |u'(x,t)|^2\,\de x+\int_0^L W(x,u(x,t))\,\de x \\
& -\frac12\int_0^L |u'(x,t-h)|^2\,\de x-\int_0^L W(x,u(x,t-h))\,\de x \\
={} & \frac12 \int_0^L \bigl( u'(x,t)-u'(x,t-h)\bigr) \bigl(u'(x,t)+u'(x,t-h)\bigr)\,\de x \\
& +\int_0^L \bigg( \int_0^1 \partial_s W(x,u_h^\lambda(x,t))\,\de\lambda\bigg)\big(u(x,t)-u(x,t-h)\big)\,\de x, \\
\intertext{which, by adding and subtracting $u'(x,t)$ in the first integral and $\partial_s W(x,u(x,t))$ in the second, becomes}
\leq{} & \int_0^L \bigl(u'(x,t)-u'(x,t-h)\bigr)u'(x,t)\,\de x \\
& +\int_0^L \partial_s W(x,u(x,t))\bigl(u(x,t)-u(x,t-h)\bigr)\,\de x + R_h^{(1)}(t), \\
\intertext{where $R_h^{(1)}(t) \coloneqq \int_0^L \bigl[ \int_0^1 \bigl( \partial_s W(x,u_h^\lambda(x,t)) - \partial_s W(x,u(x,t)) \bigr) \,\de\lambda \bigr] \bigl(u(x,t)-u(x,t-h)\bigr)\,\de x$. Integrating by parts in the first integral then yields}
={}& \int_0^L \big(u(x,t)-u(x,t-h)\big)\big(-u''(x,t)+\partial_s W(x,u(x,t))\big)\,\de x + R_h^{(1)}(t) \\
={}& \scp{u(t)-u(t-h)}{\mu(t)}_{V',V}+ R_h^{(1)}(t),
\end{split}
\end{equation*}
so that we have
\begin{subequations}\label{ee00}
\begin{equation}\label{ee01}
\mathcal{E}(u(t))-\mathcal{E}(u(t-h)) \leq \scp{u(t)-u(t-h)}{\mu(t)}_{V',V}+ R_h^{(1)}(t).
\end{equation}
By adding and subtracting $u'(x,t-h)$ and $\partial_s W(x,u(x,t-h))$ instead, an analogous computation gives
\begin{equation}\label{ee02}
\mathcal{E}(u(t))-\mathcal{E}(u(t-h)) \geq \scp{u(t)-u(t-h)}{\mu(t-h)}_{V',V}+ R_h^{(2)}(t),
\end{equation}
\end{subequations}
where $R_h^{(2)}(t) \coloneqq \int_0^L \bigl[ \int_0^1 \bigl( \partial_s W(x,u_h^\lambda(x,t)) - \partial_s W(x,u(x,t-h)) \bigr) \,\de\lambda \bigr] \bigl(u(x,t)-u(x,t-h)\bigr)\,\de x$. We now reproduce an argument similar to \cite[Proposition~4.2]{CKRS07}. 
From \eqref{ee00} we conclude \eqref{nd362} if we prove that, taking the limit $h\to0$,
\begin{subequations}\label{nd371bis}
\begin{align}
\frac{1}{h}\int_h^T \scp{u(t)-u(t-h)}{\mu(t-h)}_{V',V} v(t)\,\de t
& = \frac{1}{h}\int_0^{T-h} \scp{u(t+h)-u(t)}{\mu(t)}_{V',V} v(t+h)\,\de t \nonumber \\
& \to \int_0^T \scp{\dot{u}(t)}{\mu(t)}_{V',V} v(t)\,\de t \,, \label{nd371abis}\\
\frac{1}{h}\int_h^T \scp{u(t)-u(t-h)}{\mu(t)}_{V',V} v(t)\,\de t
& \to \int_0^T \scp{\dot{u}(t)}{\mu(t)}_{V',V} v(t)\,\de t \,, \label{nd371bbis}\\
\frac{1}{h}\int_h^T \bigl( \mathcal{E}(u(t))-\mathcal{E}(u(t-h)) \bigr) v(t)\,\de t
& = \frac{1}{h}\int_0^T \mathcal{E}(u(t)) \bigl( v(t) - v(t+h) \bigr) \,\de t \nonumber \\
& \to - \int_0^T \mathcal{E}(u(t)) \dot{v}(t)\,\de t\,, \label{nd371cbis} \\
\frac{1}{h}\int_h^T R_h^{(i)}(t)v(t)\,\de t &\to0, \qquad i=1,2, \label{ee03}
\end{align}
\end{subequations}
where $v\in W^{1,\infty}(\mathbb{R})$ is any nonnegative Lipschitz function such that $\supp v\subset[h,T-h]$.
To prove the first convergence \eqref{nd371abis}, we observe that by \cite[Chapter~1, Theorem~3.2]{B}
\begin{equation}\label{ee05}
\frac{u(t+h)-u(t)}{h} \to \dot{u}(t) \qquad\text{in $V'$, for a.e. }t\in(0,T),
\end{equation}
and therefore
\begin{equation*}
\frac1h \scp{u(t+h)-u(t)}{\mu(t)}_{V',V} v(t+h) \chi_{(0,T-h)}(t) \to \scp{\dot{u}(t)}{\mu(t)}_{V',V} v(t)
\end{equation*}
almost everywhere in $(0,T)$. Furthermore, for every measurable $A\subset[0,T]$ we have
\begin{align*}
\int_A \bigg| \frac1h \langle u(t+h)&-u(t),\mu(t) \rangle_{V',V} \,v (t+h) \chi_{(0,T-h)}(t) \bigg|  \,\de t \\
& \leq \|v\|_{L^\infty(0,T)} \biggl(\int_0^{T-h} \frac{1}{h^2}\|u(t+h)-u(t)\|^2_{V'}\,\de t\biggr)^\frac12 \biggl(\int_A \|\mu(t)\|^2_V \,\de t\biggr)^{\frac12} \\
& \leq C \|\mu\|_{L^2(A;V)}\,,
\end{align*}
where the integral of the difference quotient of $u$ is uniformly bounded by \cite[Chapter~1, Theorem~3.3]{B}. The previous estimate implies equi-integrability and therefore convergence \eqref{nd371abis} follows by Vitali's Convergence Theorem. 
The proof of \eqref{nd371bbis} is analogous. For \eqref{nd371cbis} we only need that $\mathcal{E}(u(t))\in L^1(0,T)$; in turn, this follows easily from the condition $u\in C([0,T];V)$, see Remark~\ref{rmk:interpolation}. Finally, to show \eqref{ee03}, we first observe that by the regularity \eqref{nd121} of the potential $W$ and the uniform bound
\begin{equation*}
\sup_{t\in(0,T)} \|u(t)\|_{L^\infty(0,L)} <\infty
\end{equation*}
(which follows also from Remark~\ref{rmk:interpolation}), we can find a constant $C>0$, independent of $t\in[h,T]$ and $x\in[0,L]$, such that
\begin{equation*}
|\partial_s W(x,u_h^\lambda(x,t))-\partial_s W(x,u(x,t))| \leq C(1-\lambda)|u(x,t)-u(x,t-h)|\,.
\end{equation*}
Therefore we can estimate
\begin{equation*}
\begin{split}
\frac{1}{h}\int_h^T R_h^{(1)}(t)v(t)\,\de t
& \leq \frac{C}{h} \|v\|_{\infty}\int_h^T \int_0^L |u(x,t)-u(x,t-h)|^2\,\de x \,\de t \\
& \leq  C \|v\|_{\infty} \biggl(\int_h^{T} \bigg\| \frac{u(t)-u(t-h)}{h}\bigg\|_{V'}^2\,\de t \biggr)^\frac12 \biggl(\int_h^T \| u(t)-u(t-h)\|^2_V \,\de t\biggr)^{\frac12}.
\end{split}
\end{equation*}
The first integral on the right-hand side is uniformly bounded, in view of \cite[Chapter~1, Theorem~3.3]{B}; furthermore, as $u$ is continuous as a map from $[0,T]$ with values in $V$ (by Remark~\ref{rmk:interpolation}), the second integral tends to zero as $h\to0$. This proves \eqref{ee03}.

In view of \eqref{ee00}, the convergences in \eqref{nd371bis} imply that
\begin{equation*}
- \int_0^T \mathcal{E}(u(t)) \dot{v}(t)\,\de t = \int_0^T\scp{\dot{u}(t)}{\mu(t)}_{V',V} v(t)\,\de t
\end{equation*}
for every nonnegative $v\in W^{1,\infty}(\mathbb{R})$ with compact support in $(0,T)$; since both the positive and negative parts of a Lipschitz function are Lipschitz continuous, we obtain \eqref{nd362}.

\begin{remark}[Improved estimates]\label{rmk:bounds}
We now show how to obtain uniform bounds in the space $V$ for every weak solution $u(t)$, depending on the $H^1$-norm of the corresponding initial datum $u_0$. 
We can estimate the right-hand side of the energy equality \eqref{nd360} by means of Poincar\'e and Young inequalities, to obtain that
\begin{equation*}
\begin{split}
\frac{\de}{\de t}\mathcal{E}(u(t)) + \|\mu(t)\|^2_V
& = -\beta\int_0^L u'(x,t)\mu(x,t)\,\de x
\leq \beta L\norma{u(t)}_V\norma{\mu(t)}_V \\
& \leq \frac{\beta^2L^2}{2}\norma{u(t)}_V^2 + \frac12\norma{\mu(t)}_V^2 \\
& = \beta^2L^2\mathcal{E}(u(t)) -\beta^2L^2\int_0^L W(x,u(x,t))\,\de x + \frac12\norma{\mu(t)}_V^2,
\end{split}
\end{equation*}
that is, in view of \eqref{nd122},
\begin{equation*}
\frac{\de}{\de t}\mathcal{E}(u(t)) + \frac12\|\mu(t)\|^2_V \leq  \beta^2L^2\mathcal{E}(u(t)) + \beta^2L^3K_0.
\end{equation*}
Hence Gr\"{o}nwall's Lemma yields for all $t\in[0,T]$
\begin{equation*}
\mathcal{E}(u(t)) \leq \bigl( \mathcal{E}(u_0) + L K_0 \bigr) e^{\beta^2L^2T} - L K_0 .
\end{equation*}
In particular, using \eqref{nd122} once more, we conclude that
\begin{equation} \label{S22}
\norma{u(t)}_V^2 \leq 2\bigl( \mathcal{E}(u_0) + L K_0 \bigr) e^{\beta^2L^2T}  .
\end{equation}
\end{remark}

\subsection{Continuous dependence on initial data}\label{sect:contdep}
We conclude the proof of Theorem~\ref{thm:existence} by showing that estimate \eqref{S1} holds. 
We fix two initial data $u_0,\bar u_0\in V$ with $\norma{u_0}_V,\norma{\bar u_0}_V\leq M$, where $M$ is a fixed constant, and consider the difference $w(t)\coloneqq u(t)-\bar u(t)$ of the corresponding weak solutions. We first observe that, denoting by $\mu$ and $\bar{\mu}$ the functions associated to $u$ and $\bar{u}$ according to \eqref{nd356}, respectively, and by $\mu_w\coloneqq \mu-\bar{\mu}$ their difference, we have that
\begin{equation}\label{S24}
-w''(x,t) = \mu_w(x,t) - \partial_sW(x,u(x,t)) + \partial_sW(x,\bar{u}(x,t)) ,
\end{equation}
and, since by the properties of weak solutions all the functions $u(t)$, $\bar{u}(t)$, $\mu(t)$, $\bar{\mu}(t)$ belong to the space $V$, it is easily seen that the right-hand side of \eqref{S24} vanishes at $x=0$. 
By the regularity properties of weak solutions, this implies that $w''(t)\in V$ and is therefore an admissible test function in the weak equation \eqref{nd358}. As both functions $u$ and $\bar{u}$ satisfy \eqref{nd358}, by subtracting the two resulting equations tested with $-w''(t)$ we have
\begin{equation}\label{S25}
-\scp{\dot w(t)}{w''(t)}_{V',V}
= (\mu_w(t),w''(t))_V + \beta\int_0^L w'(x,t)w''(x,t)\,\de x \,.
\end{equation}
By introducing for notational convenience the function
\begin{equation*}
\nu(x,t) \coloneqq \partial_sW(x,u(x,t)) - \partial_sW(x,\bar{u}(x,t)) ,
\end{equation*}
equation \eqref{S25} gives
\begin{equation} \label{S27}
\begin{split}
-\scp{\dot w(t)}{w''(t)}_{V',V}
={} &  - \norma{\mu_w(t)}_V^2 + ( \mu_w(t),\nu(t) )_V - \beta\int_0^L w'\mu_w\,\de x + \beta\int_0^Lw'\nu\,\de x \\
\leq{} &  - \norma{\mu_w(t)}_V^2 + \frac12\norma{\mu_w(t)}_V^2 + \frac12\norma{\nu(t)}_V^2 +\frac12\norma{\mu_w(t)}_V^2 \\
&  + \frac{\beta^2L^2}{2}\norma{w(t)}_V^2 + \beta\norma{w(t)}_V\norma{\nu(t)}_{L^2(0,L)} .
\end{split}
\end{equation}
The left-hand side of \eqref{S27} is
\begin{equation} \label{S28}
-\scp{\dot w(t)}{w''(t)}_{V',V} = \frac12\frac{\de}{\de t}\norma{w(t)}_V^2 \,,
\end{equation}
which can be proved in the same way as \eqref{nd362}, using convergence \eqref{ee05} to obtain the limit of the incremental quotients.
We now estimate the terms on the right-hand side of \eqref{S27} containing the function $\nu(t)$. Notice first that, as $u,\bar{u}\in L^\infty(0,T;V)$, we can bound
\begin{equation} \label{S32}
|\nu(x,t)| = \big| \partial_sW(x,u(x,t)) - \partial_sW(x,\bar{u}(x,t)) \big| \leq C |w(x,t)| ,
\end{equation}
and similarly
\begin{equation} \label{S33}
\begin{split}
|\nu'(x,t)| &\leq \big| \partial^2_{sx}W(x,u(x,t)) - \partial^2_{sx}W(x,\bar{u}(x,t)) \big| \\
& \qquad + \big| \partial^2_{ss}W(x,u(x,t)) \big| |u'(x,t)-\bar{u}'(x,t)| \\
& \qquad + \big| \partial^2_{ss}W(x,u(x,t))-\partial^2_{ss}W(x,\bar{u}(x,t))\big| |\bar{u}'(x,t)| \\
& \leq C |w(x,t)| + C|w'(x,t)| + C|w(x,t)||\bar{u}'(x,t)| .
\end{split}
\end{equation}
The constant $C$ appearing in \eqref{S32} and \eqref{S33} is bounded by the $C^3$-norm of the potential $W$ on the compact set
\begin{equation*}
[0,L] \times \Bigl\{s\in\R{} : |s|\leq \sup_{t\in(0,T)}\max{\{ \norma{u(t)}_{L^\infty(0,L)}, \norma{\bar{u}(t)}_{L^\infty(0,L)} \}}\Bigr\},
\end{equation*}
and therefore, in view of estimate \eqref{S22}, depends ultimately only on the fixed parameters of the problem ($L$, $\beta$, $T$, $W$) and on the energy of the initial data $\mathcal{E}(u_0)$, $\mathcal{E}(\bar{u}_0)$, which is in turn bounded uniformly by the constant $M$.

Then for the $L^2$-norm of $\nu(t)$ we have, by \eqref{S32},
\begin{equation} \label{S29}
\norma{\nu(t)}_{L^2(0,L)} \leq C\norma{w(t)}_{L^2(0,L)} ,
\end{equation}
while for the norm of $\nu(t)$ in $V$ we have, by \eqref{S33},
\begin{equation} \label{S30}
\norma{\nu(t)}_{V} \leq C\norma{w(t)}_{L^2(0,L)} +  C\norma{w(t)}_{V} + C\norma{w(t)}_{L^\infty(0,L)}\norma{\bar{u}(t)}_{V} \leq C\norma{w(t)}_V
\end{equation}
(where the norm of $\bar{u}(t)$ in $V$ can be bounded using \eqref{S22} by a constant with the same properties as before). By inserting \eqref{S28}, \eqref{S29}, and \eqref{S30} into \eqref{S27} we end up with
\begin{equation*}
\frac12\frac{\de}{\de t}\norma{w(t)}_V^2
\leq \frac12\norma{\nu(t)}_V^2 + \frac{\beta^2L^2}{2}\norma{w(t)}_V^2 + \beta\norma{w(t)}_V\norma{\nu(t)}_{L^2(0,L)} \leq C\norma{w(t)}_V^2,
\end{equation*}
and Gr\"onwall's Lemma yields
\begin{equation*}
\norma{w(t)}_V^2 \leq \norma{w(0)}_V^2 e^{CT}.
\end{equation*}
As observed before, the constant $C$ appearing in the previous estimate depends only on $L$, $\beta$, $T$, $W$, and on the uniform bound $M$ on the $V$-norm of $u_0$ and $\bar{u}_0$. This remark completes the proof of the continuous dependence of the solution on the initial datum.

\section{The case \texorpdfstring{$\beta$}{beta} small} \label{sect:smallbeta}
In this section, we first show in Theorem~\ref{thm:gb-gen-semiflow} the existence of a global attractor for weak solutions $u$ to \eqref{sub2}. In Theorem~\ref{prop:beta} we prove that, as $\beta\to 0$, solutions to \eqref{sub2} converge to that of a Cahn-Hilliard-type equation, with a quantitative estimate of the convergence rate in terms of $\beta$.

\subsection{Preliminary results on global attractors} \label{subs:prel-gen-sem}
We first recall for the reader's convenience some basic notation and results from the classical theory of attractors for semigroups (see \cite{Hale, Lad}). In what follows, let $(X,d)$ be a complete metric space. 

\begin{defin}[Semigroup]\label{def:semigroup}
A family $\mathscr{G}=\{ g_t : t\in[0,\infty)\}$ of maps $g_t:X\to X$ is a (continuous) \textit{semigroup} if the following conditions are satisfied:
\begin{itemize}
\item[(sg1)] $g_{t_1}(g_{t_2}(x))=g_{t_1+t_2}(x)$ for all $t_1,t_2\geq0$ and $x\in X$;
\item[(sg2)] $g_0(x)=x$ for all $x\in X$;
\item[(sg3)] the mapping $(t,x)\mapsto g_t(x)$ from $[0,\infty)\times X\to X$ is continuous.
\end{itemize}
\end{defin}

The following definition gathers some basic properties of semigroups.

\begin{defin}[Properties of semigroups] \label{def:prop-semigroup}
Let $\mathscr{G}=\{ g_t : t\in[0,\infty)\}$ be a semigroup.

$\mathscr{G}$ is \textit{bounded} if the set
$$
\bigcup_{x\in B}\bigl\{ g_t(x) \,:\, t\geq0 \bigr\}
$$
is bounded for every bounded set $B\subset X$.

$\mathscr{G}$ is \textit{point dissipative} if there exists a bounded set $B_0\subset X$ which attracts each point of the space, i.e., for every $x\in X$ there holds $g_t(x)\in B_0$ for every $t\geq t_x$, for some $t_x\geq 0$.

$\mathscr{G}$ is \textit{compact} if for every $t>0$ the operator $g_t$ is compact, i.e., for every bounded set $B\subset X$ its image $g_t(B)$ is precompact.
\end{defin}

A set $E\subset X$ is \textit{invariant} for the semigroup $\mathscr{G}$ if $g_t(E)=E$ for all $t>0$. The {\it $\omega$-limit} of a set $E\subset X$ is defined as the set
\begin{align*}
\omega(E)\coloneqq \bigl\{ w\in X:\,\exists \{x_k\}\subset E\text{ and } t_k\to\infty \text{ with } g_{t_k}(x_k)\to w \bigr\}.
\end{align*}
We finally recall the notion of global attractor.

\begin{defin}[Global attractor]
A set $\mathscr{A}\subset X$ is the \textit{global attractor} for a semigroup $\mathscr{G}$ if it is the minimal closed set which attracts all bounded subsets of $X$, namely for every bounded $B\subset X$ there holds $\mathrm{dist}(g_t(B);\mathscr{A})\to 0$ as $t\to +\infty$.
\end{defin}

The following classical result provides a criterion for the existence of a global attractor (see for instance \cite{Lad}).

\begin{theorem}[Existence of the global attractor] \label{thm:global-attractor}
Suppose that the semigroup $\mathscr{G}$ is bounded, point dissipative, and compact.
Then $\mathscr{G}$ has a global attractor $\mathscr{A}$, which is compact and invariant. The global attractor $\mathscr{A}$ is unique, and given by the set
$$\mathscr{A}\coloneqq\bigcup\{\omega(B):\,B\text{ is a bounded subset of }\,X\}.$$
Additionally, $\mathscr{A}$ is the maximal compact invariant subset of $X$.
\end{theorem}

\subsection{Existence of the global attractor}
We proceed by adapting the formalism introduced in Subsection \ref{subs:prel-gen-sem} to our setting.
Thanks to Theorem~\ref{thm:existence}, for $\beta\geq 0$ the equation \eqref{sub2} defines a continuous semigroup $\mathscr{G}_{\beta}$ on the space $V$.
Indeed, for every initial datum $u_0\in V$ we proved in Theorem~\ref{thm:existence} existence and uniqueness of a weak solution in $[0,T]$ starting from $u_0$, for every $T>0$; by the arbitrariness of $T>0$, this provides existence and uniqueness of a continuous map $u\colon [0,+\infty)\to V$, $u(0)=u_0$, satisfying the following properties:
\begin{itemize}
\item[(sgb1)] $u\in H^1_\textrm{loc}(0,+\infty;V')\cap L^2_\textrm{loc}(0,+\infty;H^3(0,L))$;
\item[(sgb2)] $\mu= -u'' + \partial_sW(x,u)\in L^2_\textrm{loc}(0,+\infty;V)$;
\item[(sgb3)] for every $\psi\in V$ and almost every $t>0$ we have 
$$\scp{\dot u(t)}{\psi}_{V',V}= - (\mu(t),\psi)_V-\beta\int_0^L u'(x,t)\psi(x)\,\de x;$$
\item[(sgb4)] for almost every $t>0$, $u$ satisfies the boundary conditions $u(0,t)=0$, $u'(L,t)=0$.
\end{itemize}
Then it is easily seen that the family $\mathscr{G}_\beta=\{ g_t : t\in[0,\infty)\}$, where, for every $u_0\in V$ and $t\geq0$, $g_t(u_0)$ is the value $u(t)$ of the unique weak solution starting from $u_0$, is a continuous semigroup. In particular, condition (sg3) in Definition~\ref{def:semigroup} is a direct consequence of \eqref{S1}. 

We are now in a position to prove the existence of a global attractor.

\begin{theorem}\label{thm:gb-gen-semiflow}
For every $\beta<L^{-3}$, the semigroup $\mathscr{G}_{\beta}$ has a global attractor.
\end{theorem}

\begin{proof}
The proof consists in showing that for $\beta$ small enough the semigroup $\mathscr{G}_{\beta}$ satisfies the conditions of Theorem~\ref{thm:global-attractor} characterizing the existence of a global attractor. For convenience of the reader, we subdivide the proof into two steps.

\smallskip
\noindent\textit{Step 1: $\mathscr{G}_{\beta}$ is bounded and point dissipative.}
Let $u_0\in V$, let $u$ be the weak solution to \eqref{sub2} with initial datum $u_0$ provided by Theorem~\ref{thm:existence}, and fix $t\in [0,+\infty)$. We first observe that, by \eqref{nd122} and using Young's and Poincar\'e's inequalities, we can estimate
\begin{equation} \label{eq:1ga}
\begin{split}
\mathcal{E}(u(t))
& \leq - \int_0^L u''(x,t)u(x,t)\,\de x + \int_0^L \partial_s W(x,u(x,t))u(x,t)\,\de x + K_1 L - \frac12\|u(t)\|_V^2 \\
& = \int_0^L \mu(x,t) u(x,t)\, \de x + K_1L - \frac12\|u(t)\|_V^2 \\
& \leq L^2\|\mu(t)\|_{V} \|u(t)\|_{V} + K_1L - \frac12\|u(t)\|_V^2
\leq \frac{L^4}2 \|\mu(t)\|_V^2 + K_1 L \,.
\end{split}
\end{equation}
By plugging this estimate into the energy equality \eqref{nd360} we obtain 
\begin{equation*}
\begin{split}
\frac{\de}{\de t}\mathcal{E}(u(t)) + \frac{1}{L^4} \mathcal{E}(u(t)) + \frac12\|\mu(t)\|_V^2
& \stackrel{\eqref{eq:1ga}}{\leq} \frac{\de}{\de t}\mathcal{E}(u(t)) +\|\mu(t)\|_V^2 + \frac{K_1}{L^3} \\
& \stackrel{\eqref{nd360}}{=} -\beta\int_0^L u'(x,t)\mu(x,t)\,\de x +  \frac{K_1}{L^3} \\
& \leq \frac{\beta^2 L^2}{2} \|u(t)\|_V^2 + \frac12\|\mu(t)\|_V^2 + \frac{K_1}{L^3} \\
& \stackrel{\eqref{nd122}}{\leq} \beta^2 L^2\mathcal{E}(u(t)) + \frac12\|\mu(t)\|_V^2 + \frac{K_1}{L^3} + \beta^2 L^3 K_0 \,,
\end{split}
\end{equation*}
which yields
\begin{equation*}
\frac{\de}{\de t}\mathcal{E}(u(t)) \leq -C_\beta \mathcal{E}(u(t)) + D_\beta,
\end{equation*}
with the positions
\begin{equation*}
C_\beta\coloneqq \frac{1}{L^4}-\beta^2L^2 \qquad\text{and}\qquad D_\beta\coloneqq \frac{K_1}{L^3} + \beta^2 L^3 K_0 .
\end{equation*}
Notice that $C_\beta>0$ for $\beta < \frac{1}{L^3}$; therefore Gr\"onwall's inequality implies
\begin{equation} \label{eq:5ga}
\mathcal{E}(u(t))\leq \Big(\mathcal{E}(u_0)-\frac{D_\beta}{C_\beta}\Big)e^{-C_\beta t}+\frac{D_\beta}{C_\beta},
\end{equation}
for every $t\in [0,+\infty)$. The fact that $\mathscr{G}_{\beta}$ is point dissipative and bounded follows then by \eqref{eq:5ga} and \eqref{nd359}.

\smallskip
\noindent\textit{Step 2: $\mathscr{G}_{\beta}$ is compact.}
Let $B\subset V$ be a bounded set of initial data, and consider the family $S\coloneqq\{ u(\cdot): u(0)\in B \}$ of the weak solutions starting from a point in $B$. We first remark that the set $S$ is bounded in $L^2(0,T;H^3(0,L))\cap H^1(0,T;V')$, for every $T>0$. This follows from the construction of the weak solution that was performed in Sections~\ref{sect:MM} and~\ref{sect:fixedpoint}: indeed, every weak solution $u$ fulfills by construction estimates \eqref{nd190}--\eqref{nd191}, that is,
\begin{equation*}
\|u\|_{L^2(0,T;H^3(0,L))} + \|\dot{u}\|_{L^2(0,T;V')} \leq M
\end{equation*}
where $M$ depends on $\beta$, $L$, $T$, $W$, $\mathcal{E}(u(0))$, and $\|u\|_{L^2(0,T;V)}$. As these quantities are uniformly bounded with respect to $u\in S$ thanks to the assumption that $\|u(0)\|_V$ are uniformly bounded (recall also \eqref{S22}), we can conclude that the set $S$ is bounded in $L^2(0,T;H^3(0,L))\cap H^1(0,T;V')$, for every $T>0$.

By Aubin-Lions-Simon Theorem (see \cite[Theorem~8.62]{Leo}), for every $T>0$ the space $L^2(0,T;H^3(0,L))\cap H^1(0,T;V')$ is compactly embedded in $L^2(0,T;H^2(0,L))$, and therefore by a diagonal argument we can extract a subsequence $\{u_j\}_{j}\subset S$ such that
\begin{equation*}
u_{j} \to u \qquad\text{in }L^2(0,T;H^2(0,L)), \text{ for all $T>0$,}
\end{equation*}
as $j\to\infty$, for some $u\in L^2_{\mathrm{loc}}(0,\infty;H^2(0,L))$. In view of Remark~\ref{rmk:interpolation}, $u\in C(0,\infty;V)$. Up to further subsequences, we can also assume that
\begin{equation} \label{eq:1ac}
u_{j}(t) \to u(t) \qquad\text{in $H^2(0,L)$, for a.e. $t>0$.}
\end{equation}
The goal is now to show that the convergence $u_{j}(t)\to u(t)$ in $V$ holds \textit{for every} $t>0$.

We first claim that
\begin{equation} \label{eq:6ac}
u_{j}(t) \wto u(t) \qquad\text{weakly in $V$, for every $t>0$.}
\end{equation}
Indeed, for any given $T>0$ and for all $0\leq s \leq t \leq T$ we have
\begin{equation*}
\begin{split}
\sup_j \norma{u_{j}(t)-u_{j}(s)}_{V'}^2
& \leq \sup_j \biggl(\int_s^t \norma{\dot u_{j}(\tau)}_{V'}\,\de\tau\biggr)^2 \\
& \leq (t-s)\sup_j \int_s^t\norma{\dot u_{j}(\tau)}_{V'}^2\de\tau
\leq (t-s)\sup_j\|u_{j}\|_{H^1(0,T;V')}^2,
\end{split}
\end{equation*}
which shows that the sequence of functions $\{u_{j}\}_j$ is uniformly equicontinuous from $[0,T]$ to $V'$. In addition, by \eqref{S22} the functions $u_{j}(t)$ are equibounded in $V$, and therefore $\{u_{j}(t)\}_j$ is relatively compact in $V'$ for all $t\in[0,T]$. Hence by Ascoli-Arzel\`a Theorem we can extract a further subsequence such that
\begin{equation} \label{eq:8ac}
u_{j_{k}}(t) \to v(t) \qquad\text{in $V'$, for every $t\in[0,T]$,}
\end{equation}
for some $v\in C(0,T;V')$. In view of \eqref{eq:1ac} and the fact that $u\in C(0,\infty;V)$, we necessarily have $u=v$ and the convergence \eqref{eq:8ac} holds for the full sequence $(u_{j})_j$. Finally, for every $t>0$ the sequence $(u_{j}(t))_j$ is bounded in $V$, and therefore by \eqref{eq:8ac} it converges weakly in $V$ to $u(t)$. Then the claim \eqref{eq:6ac} is proved.

The final step amounts to show that the convergence in \eqref{eq:6ac} is actually strong in the space $V$. To this aim, we introduce the functions $\mu_{j}$ and $\mu$ associated to $u_{j}$ and $u$ by \eqref{nd356}, respectively, and we notice that, by \eqref{eq:1ac} and the continuity of $\partial_s W$, we have
\begin{equation} \label{eq:4ac}
\mu_{j}(t)\to\mu(t) \qquad\text{in $L^2(0,L)$, for a.e. $t>0$.}
\end{equation}
In particular, the convergences \eqref{eq:1ac} and \eqref{eq:4ac} guarantee that for every $t>0$
\begin{equation} \label{eq:5ac}
\lim_{j\to\infty} \int_0^t\int_0^L u_{j}'(x,\tau)\mu_{j}(x,\tau)\,\de x \de \tau = \int_0^t\int_0^L u'(x,\tau)\mu(x,\tau)\,\de x \de \tau\,.
\end{equation}
We now let
\begin{equation*}
\widetilde{\mathcal{E}}(u(t)) \coloneqq \mathcal{E}(u(t)) + \beta\int_0^t\int_0^L u'(x,\tau)\mu(x,\tau)\,\de x \de \tau\,,
\end{equation*}
where $\mathcal{E}$ is defined in \eqref{nd359}. By the energy equality \eqref{nd360}, the map $t\mapsto\widetilde{\mathcal{E}}(u_{j}(t))$ is monotone non-increasing for every $j$. Moreover, by \eqref{eq:1ac} and \eqref{eq:5ac} there holds $\widetilde{\mathcal{E}}(u_{j}(t))\to\widetilde{\mathcal{E}}(u(t))$ for almost every $t>0$. Therefore by monotonicity and continuity of $\widetilde{\mathcal{E}}(u(\cdot))$ we have
\begin{equation*}
\widetilde{\mathcal{E}}(u_{j}(t))\to\widetilde{\mathcal{E}}(u(t)) \qquad\text{for every }t>0.
\end{equation*}
By using again \eqref{eq:5ac} this yields $\mathcal{E}(u_{j}(t))\to\mathcal{E}(u(t))$ for every $t>0$ and, in turn,
\begin{equation} \label{eq:7ac}
\begin{split}
\lim_{j\to\infty} \biggl( \frac12\|u_{j}(t)\|_{V}^2 - \frac12\|u(t)\|_V^2 \biggr)
& = \lim_{j\to\infty} \int_0^L \biggl( W(x,u(x,t))-W(x,u_{j}(x,t))\biggr)\,\de x
\end{split}
\end{equation}
for every $t>0$. Since the weak convergence in $V$ implies uniform convergence in $[0,L]$, by continuity of the potential $W$ it is easily seen that the condition \eqref{eq:6ac} implies that the limit on the right-hand side of \eqref{eq:7ac} is zero for every $t>0$, and therefore $\|u_{j}(t)\|_V\to\|u(t)\|_V$ for every $t>0$. This information, combined with \eqref{eq:6ac}, allows us to conclude that $u_{j}(t)\to u(t)$ strongly in $V$ for every $t>0$. This concludes the proof of the compactness of $\mathscr{G}_\beta$.
\end{proof}

\subsection{Further characterization of solutions}
In the limit $\beta\to0$ we can show that weak solutions to \eqref{sub2} are close, in the norm of the space $V$, to the weak solution of the limit problem with $\beta=0$, corresponding to the same initial datum. More precisely, the following result holds true.

\begin{theorem}\label{prop:beta}
Let $u_0\in V$ be a fixed initial datum, and let $T>0$, $\beta_0>0$ be given. For $\beta\in(0,\beta_0)$, we denote by $u_\beta$ the unique weak solution to \eqref{sub2} in $[0,T]$ corresponding to $u_0$, constructed in Theorem~\ref{thm:existence}. Moreover, we denote by $\bar{u}$ the unique weak solution with the same initial datum $u_0$ corresponding to $\beta=0$. Then there exists a constant $C>0$ (depending on $L$, $W$, $T$, $\beta_0$, and on $\|u_0\|_V$) such that for all $\beta\in(0,\beta_0)$ and $t\in[0,T]$
\begin{equation} \label{bs001}
\norma{u_{\beta}(t)-\bar{u}(t)}_{V}\leq C\beta.
\end{equation}
\end{theorem}

\begin{proof}
We proceed similarly to the proof of the continuous dependence of solutions on the initial data, see Subsection~\ref{sect:contdep}.
We first observe that the difference of the two solutions $w\coloneqq u_\beta-\bar{u}$ solves
\begin{equation}\label{bs020}
\scp{\dot w(t)}{\psi}_{V',V}
= -(\mu_w(t),\psi)_V - \beta\int_0^L u_\beta'(x,t)\psi(x)\,\de x
\end{equation}
for every test function $\psi\in V$, where we denoted by $\mu_w\coloneqq \mu_\beta-\bar{\mu}$ the difference of the maps $\mu_\beta$ and $\bar{\mu}$ associated to $u_\beta$ and $\bar{u}$ according to \eqref{nd356}, respectively. Notice that
\begin{equation}\label{bs021}
-w''(x,t) = \mu_w(x,t) - \nu(x,t), 
\end{equation}
where we defined $\nu(x,t) \coloneqq \partial_sW(x,u_\beta(x,t)) - \partial_sW(x,\bar{u}(x,t))$. From \eqref{bs021} it is easily seen that $w''(t)\in V$ and is therefore an admissible test function in the weak equation \eqref{bs020}:
\begin{equation}\label{bs022}
-\scp{\dot w(t)}{w''(t)}_{V',V}
= (\mu_w(t),w''(t))_V + \beta\int_0^L u_\beta'(x,t)w''(x,t)\,\de x \,.
\end{equation}
By arguing as in \eqref{S27}--\eqref{S28}, equation \eqref{bs022} gives
\begin{equation} \label{bs024}
\begin{split}
\frac12\frac{\de}{\de t}\norma{w(t)}_V^2
& = - \norma{\mu_w(t)}_V^2 + ( \mu_w(t),\nu(t) )_V - \beta\int_0^L u_\beta'\mu_w\,\de x + \beta\int_0^Lu_\beta'\nu\,\de x \\
& \leq - \norma{\mu_w(t)}_V^2 + \frac12\norma{\mu_w(t)}_V^2 + \frac12\norma{\nu(t)}_V^2 +\frac12\norma{\mu_w(t)}_V^2 \\
& \qquad\qquad\qquad\qquad + \frac{\beta^2L^2}{2}\norma{u_\beta(t)}_V^2 + \beta\norma{u_\beta(t)}_V\norma{\nu(t)}_{L^2(0,L)} .
\end{split}
\end{equation}
Notice now that, thanks to \eqref{S22}, we have a uniform estimate
\begin{equation*}
\|u_\beta(t)\|_V^2 \leq C_1,
\end{equation*}
where $C_1$ depends ultimately only on $L$, $W$, $T$, $\beta_0$, and on $\|u_0\|_V$. In turn, by arguing as in \eqref{S29}--\eqref{S30} the norm of the function $\nu(t)$ can be bounded as
\begin{equation*}
\norma{\nu(t)}_{V}^2 \leq C_2\norma{w(t)}_V^2
\end{equation*}
for another constant $C_2$ depending on the same quantities as $C_1$.
Therefore \eqref{bs024} yields
\begin{equation*}
\frac12\frac{\de}{\de t}\norma{w(t)}_V^2
\leq \norma{\nu(t)}_V^2 + \beta^2L^2\norma{u_\beta(t)}_V^2
\leq C_2\norma{w(t)}_V^2 + \beta^2C_1L^2,
\end{equation*}
and, as $w(0)=0$, an application of Gr\"onwall's Lemma yields \eqref{bs001}.
\end{proof}

\appendix
\section{} \label{sect:appendix}

We collect here a few results of technical nature which are needed throughout the paper.
The following lemma will be instrumental in the proof of Proposition~\ref{prop:compact}.

\begin{lemma} \label{lem:q}
Let $q\in L^2(0,T;V)$ and let $\tau \coloneqq \frac{T}{N}$ for $N\in\mathbb{N}$. Then, as $\tau\to0$,
\begin{equation}\label{A4}
\frac{1}{\tau}\sum_{k=1}^N \chi_{((k-1)\tau,k\tau)}(t) \int_{(k-1)\tau}^{k\tau} q'(s)\,\de s
\to q'(t) \qquad\text{weakly in } L^2(0,T;L^2(0,L)).
\end{equation}
\end{lemma}

\begin{proof}
We set
\begin{equation*}
g_\tau(x,t)\coloneqq \frac{1}{\tau}\sum_{k=1}^N \chi_{((k-1)\tau,k\tau)}(t) \int_{(k-1)\tau}^{k\tau} \bigl( q'(x,t)-q'(x,s) \bigr)\,\de s.
\end{equation*}
We can estimate
\begin{align*}
|g_\tau(x,t)|^2
& = \frac{1}{\tau^2} \sum_{k=1}^N \chi_{((k-1)\tau,k\tau)}(t) \bigg| \int_{(k-1)\tau}^{k\tau} \bigl( q'(x,t)-q'(x,s) \bigr)\,\de s \bigg|^2 \\
& \leq \frac{1}{\tau} \sum_{k=1}^N \chi_{((k-1)\tau,k\tau)}(t) \int_{(k-1)\tau}^{k\tau} |q'(x,t)-q'(x,s)|^2 \,\de s \\
& \leq 2|q'(x,t)|^2 + \frac{2}{\tau} \sum_{k=1}^N \chi_{((k-1)\tau,k\tau)}(t) \int_{(k-1)\tau}^{k\tau} |q'(x,s)|^2 \,\de s\,,
\end{align*}
so that by integrating on $(0,L)\times(0,T)$ we obtain
\begin{equation*}
\int_0^T\int_0^L |g_\tau(x,t)|^2\,\de x\de t \leq 4\int_0^T\int_0^L |q'(x,t)|^2\,\de x \de t =:M <\infty.
\end{equation*}
This implies that, up to subsequences, $g_\tau\wto g$ weakly in $L^2((0,L)\times(0,T))$.

We shall now show that $g=0$. The estimate
\begin{equation*} 
|g_\tau(x,t)| \leq \frac{1}{\tau} \int_{t-\tau}^{t+\tau} \big| q'(x,t)-q'(x,s) \big| \,\de s
\end{equation*}
implies that $\|g_\tau(t)\|_{L^2(0,L)}\to0$ for almost every $t\in(0,T)$ by \cite[Theorem 8.19]{Leo}; in turn, given $\eps>0$, by Egorov's Theorem we can find a set $A_\eps\subset(0,T)$ with measure smaller than $\eps$ such that $\|g_\tau(t)\|_{L^2(0,L)}\to0$ uniformly on $(0,T)\setminus A_\eps$. Therefore for any test function $\varphi\in L^2((0,L)\times(0,T))$ we have
\begin{align*}
\limsup_{\tau\to0}\bigg|\int_0^T\int_0^L g_\tau(x,t)\varphi(x,t)\,\de x\de t \bigg|
& \leq \limsup_{\tau\to0} \int_0^T \|g_\tau(t)\|_{L^2(0,L)} \|\varphi(t)\|_{L^2(0,L)} \,\de t\\
& = \limsup_{\tau\to0} \int_{A_\eps} \|g_\tau(t)\|_{L^2(0,L)} \|\varphi(t)\|_{L^2(0,L)} \,\de t\\
& \leq M^\frac12\biggl( \int_{A_\eps} \|\varphi(t)\|_{L^2(0,L)}^2\,\de t \biggr)^\frac12,
\end{align*}
and the right-hand side can be made arbitrarily small as $\eps\to0$. This shows that  $g_\tau\wto0$ weakly in $L^2((0,L)\times(0,T))$ and \eqref{A4} holds.
\end{proof}

We discuss in the following two lemmas the dependence on $t$ of the function $z_{\varphi(t)}$ introduced in \eqref{nd108}, when $\varphi$ depends on $t$.
\begin{lemma} \label{lem:zq}
Let $q\in L^2(0,T;V')$ and let $z_{q(t)}$ be defined by \eqref{nd108}.
Then the map $t\mapsto z_{q(t)}$ belongs to $L^2(0,T;V)$.
\end{lemma}

\begin{proof}
For every $\varphi \in V'$ the map
\begin{equation*}
t \mapsto \scp{\varphi}{z_{q(t)}}_{V',V}
\stackrel{\eqref{nd108}}{=} \int_0^L z_\varphi'(x)z_{q(t)}'(x)\,\de x
= \scp{q(t)}{z_\varphi}_{V',V}
\end{equation*}
is measurable thanks to the assumption $q\in L^2(0,T;V')$. Hence the map  $t\mapsto z_{q(t)}$ is weakly measurable from $(0,T)$ to $V$ and in turn strongly measurable, by Pettis Theorem (see \cite[Theorem~8.3]{Leo}) and the separability of $V$.
Moreover, by \eqref{nd108} we have
\begin{equation*}
\int_0^T \norma{z_{q(t)}}^2_V\,\de t = \int_0^T \norma{q(t)}_{V'}^2 <\infty
\end{equation*}
since $q\in L^2(0,T;V')$, which shows that the map $t\mapsto z_{q(t)}$ is Bochner integrable and belongs to $L^2(0,T;V)$.
\end{proof}

Notice that, if $q\in L^2(0,T;L^2(0,L))$, then by \eqref{nd112bis} we have $q\in L^2(0,T;V')$. Hence we can apply the previous lemma to deduce that also in this case $z_{q(t)}\in L^2(0,T;V)$.

\begin{lemma} \label{lem:zpunto}
Let $u\in H^1(0,T;V')$ and let $z_{u(t)}$ be defined by \eqref{nd108}.
Then the map $\psi(t) \coloneqq z_{u(t)}$ is in $H^1(0,T;V)$, and
\begin{equation} \label{A2}
\dot{\psi}(t) = z_{\dot{u}(t)}.
\end{equation}
\end{lemma}
\begin{proof}
We first observe that $\psi\in L^2(0,T;V)$, thanks to Lemma~\ref{lem:zq}.
Similarly, the map $t\mapsto z_{\dot{u}(t)}$ is in $L^2(0,T;V)$, thanks to the same lemma applied to $q=\dot{u}\in L^2(0,T;V')$.

We now show equality \eqref{A2}, which will complete the proof of the lemma. By definition of weak derivative, \eqref{A2} is equivalent to show
\begin{equation} \label{A3}
\int_0^T z_{\dot{u}(t)}\varphi(t)\,\de t = - \int_0^T z_{u(t)}\dot{\varphi}(t)\,\de t
\end{equation}
for every $\varphi\in C^1_{\mathrm c}(0,T)$, where \eqref{A3} is an equality between elements of $V$.
For every $\eta\in V$ we have
\begin{equation*}
\begin{split}
\Bigl(\int_0^T z_{u(t)}\dot{\varphi}(t)\,\de t \,,\, \eta\Bigr)_V
& = \int_0^T \bigl(z_{u(t)}\dot{\varphi}(t), \eta\bigr)_V\,\de t
= \int_0^T \int_0^L z'_{u(t)}(x)\dot{\varphi}(t)\eta'(x) \,\de x \de t \\
& \stackrel{\eqref{nd108}}{=} \int_0^T\dot{\varphi}(t) \scp{u(t)}{\eta}_{V',V}\,\de t
= \scp{\int_0^T u(t)\dot{\varphi}(t)\,\de t}{\eta}_{V',V} \\
& = - \scp{\int_0^T \dot{u}(t)\varphi(t)\,\de t}{\eta}_{V',V}
= -\int_0^T \varphi(t) \scp{\dot{u}(t)}{\eta}_{V',V}\,\de t \\
& \stackrel{\eqref{nd108}}{=} -\int_0^T \int_{0}^L \varphi(t) z'_{\dot{u}(t)}(x)\eta'(x) \,\de x \de t
= -\int_0^T \bigl( \varphi(t) z_{\dot{u}(t)}, \eta \bigr)_V \de t \\
& = - \Bigl(\int_0^T \varphi(t)z_{\dot{u}(t)}\,\de t \,,\, \eta\Bigr)_V \,.
\end{split}
\end{equation*}
In the previous chain of equalities we used several times the property that the Bochner integral commutes with linear operators, see for instance \cite[Theorem~8.13]{Leo}. This proves claim \eqref{A3}.
\end{proof}

\medskip

\noindent {\bf Acknowledgments.} 
The authors thank the hospitality of the departments of mathematics of the universities of Heidelberg and Vienna, of SISSA, and of the E.~Schr\"{o}dinger Institute in Vienna, where this research was developed.
The authors are members of the GNAMPA group of INdAM. 
This work was partially supported by the 2015 GNAMPA project \emph{Fenomeni Critici nella Meccanica dei Materiali: un Approccio Variazionale}.
E.D.\@ acknowledges the support of the Austrian Science Fund (FWF) project P~27052 and of the SFB project F65 \emph{Taming complexity in partial differential systems}.
The research of M.M.\@ was partially supported by grant FCT-UTA\_CMU/MAT/0005/2009 \emph{Thin Structures, Homogenization, and Multiphase Problems}, by the ERC Advanced grant \emph{Quasistatic and Dynamic Evolution Problems in Plasticity and Fracture} (Grant agreement no.: 290888), by the ERC Starting grant \emph{High-Dimensional Sparse Optimal Control} (Grant agreement no.: 306274), and by the DFG Project \emph{Identifikation von Energien durch Beobachtung der zeitlichen Entwicklung von Systemen} (FO 767/7).
Finally, the authors thank Stefano Melchionna for fruitful discussions.

\end{document}